\documentclass[10pt]{article}

\usepackage{amsmath}
\usepackage{amsfonts}
\usepackage{amssymb}
\usepackage[english]{babel}
\usepackage{amsthm}
\usepackage{extarrows}
\usepackage{amscd}
\usepackage{tikz}
\usepackage{mathdots}
\usetikzlibrary{arrows,shapes}

\def\pf{\par\noindent {\bf Proof}~\par\noindent}
\def\qed{~\hfill{$\square$}\pagebreak[1]\par\medskip\par}

\newcommand{\mR}{\mathbb{R}}
\newcommand{\mC}{\mathbb{C}}

\newcommand{\mS}{\mathbb{S}}
\newcommand{\mH}{\mathbb{H}}
\newcommand{\mI}{\mathbb{I}}
\newcommand{\mJ}{\mathbb{J}}
\newcommand{\mK}{\mathbb{K}}
\newcommand{\mQ}{\mathbb{Q}}

\newcommand{\mcC}{\mathcal{C}}

\newcommand{\gf}{\mathfrak{f}}
\newcommand{\gfd}{\mathfrak{f}^{\dagger}}
\newcommand{\gsl}{\mathfrak{sl}}
\newcommand{\gsp}{\mathfrak{sp}}

\newcommand{\gu}{\mathfrak{u}}

\newcommand{\gosp}{\mathfrak{osp}}

\newcommand{\ol}{\overline}
\newcommand{\olz}{\ol{z}}

\newcommand{\uX}{\underline{X}}

\newcommand{\uz}{\underline{z}}
\newcommand{\uzJ}{\underline{z}^J}
\newcommand{\uzJd}{\underline{z}^{\dagger J}}
\newcommand{\uzd}{\underline{z}^{\dagger}}

\newcommand{\upz}{\partial_{\uz}}
\newcommand{\upzJ}{\partial_{\uz}^J}

\newcommand{\upzd}{\partial_{\uz}^{\dagger}}
\newcommand{\upzJd}{\partial_{\uz}^{\dagger J}}

\newcommand{\p}{\partial}
\newcommand{\dirac}{\underline{\p}}

\newcommand{\eop}{\hfill$\square$}

\newcommand{\onehalf}{\frac{1}{2}}

\newcommand{\Ker}{\rm Ker}

\newtheorem{theorem}{Theorem}
\newtheorem{lemma}{Lemma}
\newtheorem{proposition}{Proposition}
\newtheorem{definition}{Definition}
\newtheorem{remark}{Remark}
\newtheorem{corollary}{Corollary}

\begin{document}

\title{Fundaments of Quaternionic Clifford Analysis II:\\ Splitting of Equations}
\author{F.\ Brackx$^\ast$, H.\ De Schepper$^\ast$, D.\ Eelbode$^{\ast\ast}$, R.\ L\'{a}vi\v{c}ka$^\ddagger$ \& V.\ Sou\v{c}ek$^\ddagger$}

\date{\small{$^\ast$ Clifford Research Group, Dept. of Math. Analysis, 
Faculty of Engineering and Architecture,\\ Ghent University,
Building S22, Galglaan 2, B--9000 Gent, Belgium\\
$\ast\ast$ University of Antwerp, Middelheimlaan 2, Antwerpen, Belgium\\
$^\ddagger$ Charles University in Prague, Faculty of Mathematics and Physics, Mathematical Institute\\
Sokolovsk\'a 83, 186 75 Praha, Czech Republic}}

\maketitle


\begin{abstract}
Quaternionic Clifford analysis is a recent new branch of Clifford analysis, a higher dimensional function theory which refines harmonic analysis and generalizes to higher dimension the theory of holomorphic functions in the complex plane. So--called quaternionic monogenic functions satisfy a system of first order linear differential equations expressed in terms of four interrelated Dirac operators. The conceptual significance of quaternionic Clifford analysis is unraveled by 
showing that quaternionic monogenicity can be characterized by means of generalized gradients in the sense of Stein and Weiss.
At the same time, connections between quaternionic monogenic functions and other branches of Clifford analysis, viz Hermitian monogenic and standard or Euclidean monogenic functions are established as well.
\end{abstract}


\section{Introduction}
\label{introduction}


The Dirac operator plays a very important role in modern mathematics. A comprehensive description of its role in classical geometry can be found, e.g.,
in the book by Lawson and Michelsohn \cite{LM}. The Dirac operator can be characterized as the only (up to a multiple) elliptic and conformally invariant first order system of PDE's acting on spinor valued functions. A recent monograph \cite{BHMMM} on the spinorial approach to Riemannian and conformal geometry  shows the importance to considering the Dirac operator not only on Riemannian (spin) manifolds but also on manifolds with a special holonomy. By the Berger-Simons classification, we know that these include 
K\"ahler and quaternionic--K\"ahler manifolds. As shown in \cite{BHMMM}, the spectral properties of the Dirac operator on these various types of Riemannian manifolds are different.\\[-2mm]

Standard Clifford analysis, also called Euclidean or orthogonal Clifford analysis, is, in its most basic form, a higher dimensional generalization of the theory of holomorphic functions in the complex plane and at the same time a refinement of harmonic analysis, see e.g.\ \cite{bds,gilbert, dss,ghs}.  Euclidean space $\mR^m$ is the flat model for Riemannian (spin) geometry and Clifford analysis is the function theory for the Dirac operator in this flat model. We refer to this setting as the Euclidean (or orthogonal) case, since the fundamental symmetry group leaving the Dirac operator $\dirac = \sum_{\alpha=1}^m e_{\alpha}\, \p_{X_{\alpha}}$  invariant, is the special orthogonal group SO($m$),  which is doubly covered by the spin group Spin($m$) in the Clifford algebra.\\[-2mm]

The flat model for  K\"ahler  manifolds is the Euclidean space $\mR^{2n}$ of even dimension with a chosen additional datum, a  complex structure $\mI,$ i.e., an SO($2n$)--element squaring up to $-E_{2n}, E_{2n}$ being the identity matrix. There are now two invariant first order differential operators in this setting, namely, the classical Dirac operator $\dirac$ and its image under the action of the complex structure $\dirac_{\mI}=\sum_{\alpha=1}^{2n} \mI(e_{\alpha})\, \p_{X_{\alpha}},$ which are the flat versions of the invariant differential operators defined on any K\"ahler  manifold. In the books \cite{struppa,rocha} and the series of papers \cite{partI,partII, eel,eel2,sabadini} so--called Hermitian Clifford analysis emerged as a refinement of Euclidean Clifford analysis, by considering functions which now take their values in the complex Clifford algebra $\mC_{2n},$ or in the complex spinor space $\mS$. Hermitian Clifford analysis  focuses on the simultaneous null solutions of both operators $\dirac$ and $\dirac_{\mI}$, called Hermitian monogenic functions. The fundamental symmetry group underlying this function theory is the unitary group U($n$). It is worth mentioning that the traditional holomorphic functions of several complex variables appear as a special case of Hermitian monogenicity, see Section \ref{EHCA}.\\[-2mm]

The next step is to go from complex manifolds (in particular K\"ahler ones)
to manifolds with a suitable quaternionic structure. It is possible to consider in general hypercomplex manifolds, or hyper--K\"ahler manifolds, or quaternionic--K\"ahler manifolds. 
The flat model for hypercomplex manifolds is the Euclidean space $\mR^{4p}$, the dimension of which is assumed to be a fourfold, together with a choice of three complex structures $\mI,$ $\mJ$ and $\mK$ satisfying the usual multiplication relations for the three quaternionic units. It is clear that, in this setting, we have four basic operators acting on spinor valued functions, namely, $\dirac,$ $\dirac_{\mI},$ $\dirac_{\mJ}$ and $\dirac_{\mK}$. Recently a new branch of Clifford analysis arose, focussing at so-called quaternionic monogenic functions (see e.g.\ \cite{Ee,PSS,DES,paper1}).  The associated function theory is called  quaternionic Clifford analysis.
The fundamental symmetry group for quaternionic Clifford analysis is the symplectic group Sp($p$).\\[-2mm]

There is a fundamental problem with this notion of quaternionic monogenic functions, which was already visible  in the case of Hermitian Clifford analysis. In the classical orthogonal case, the Dirac operator acts naturally on spinor valued functions. The basic spinor representation is an irreducible representation of the symmetry group  SO($m$). In the Hermitian case, introducing the complex structure $\mI,$ the symmetry group is reduced to a subgroup of SO($2n$) isomorphic to U($n$). Consequently, the spinor space $\mS$ where our functions take their values is not irreducible anymore but splits into a finite sum of irreducible components instead. These components can be described as the homogeneous components $\mS^r  \simeq  \wedge^r(W^\dagger)$
of the Grassmann algebra of a maximal isotropic subspace $W^\dagger$ of the complexification $\mC^{2n}$. In the quaternionic case the reduction of the symmetry group to Sp($p$) leads to a further decomposition of the homogeneous parts $\mS^r$. The corresponding irreducible pieces $\mS^r_s$ now depend on two parameters, see Section \ref{spinordecomposition}
for a~detailed description.\\[-2mm]

A natural approach to invariant first order systems is described by Stein and Weiss in \cite{stein}, and it also works equally well on manifolds.
The structure of a hypercomplex manifold $M$ is given by a choice of  a principal U($n$)-bundle $P$ over $M$, and invariant operators in the sense of Stein and Weiss are acting on sections of bundles associated to irreducible U($n$)-modules. In our case, their scheme applies to sections  of bundles associated to the irreducible parts (with respect to Sp($p$)) of the spinor bundle $\mS$, which are induced by the representations $\mS^r_s$.  Basic first order invariant operators on hypercomplex manifolds hence can be classified in their turn by the Stein and Weiss scheme for  sections  of bundles induced by each representation $\mS^r_s$.  The corresponding first order invariant systems are the basic first order PDE's to be considered on this type of manifolds.\\[-2mm]

The basic question now is: what is the relation between quaternionic monogenic functions and the solutions of invariant first order PDE's stemming from the Stein-Weiss approach? This is the key problem in establishing quaternionic Clifford analysis, to which the present paper proposes an answer. First, in Section \ref{spinordecomposition}, we recall the decomposition of the spinor space $\mS$  into its irreducible components $\mS^r_s$. Then, in Section \ref{gradients}, we study a classification of the Stein-Weiss operators -- also called generalized gradients -- in the quaternionic setting. It leads to an alternative definition of quaternionic monogenic functions.  In Theorem \ref{main} we prove that $\mS^r_s$--valued functions are quaternionic monogenic if and only if they are common null solutions of all Stein-Weiss generalized gradients. This is the  main result of the paper and it shows that both approaches are equivalent.\\[-2mm]

In order to making the present paper self-contained, the basics of Clifford algebra are recalled in Section \ref{cliffordalgebra}.  Section \ref{EHCA} shortly outlines the basics of Euclidean and Hermitian Clifford analysis while Section \ref{QCA} introduces quaternionic Clifford analysis. In Section \ref{lowdimension} our results are illustrated by explicit calculations of the system of equations for quaternionic monogenicity in $\mR^8$. For a detailed account on the construction of the symplectic cells in spinor space, we refer to \cite{paper1}. For an account of the impact on the Fischer decomposition of harmonic polynomials of imposing symmetry with respect to the symplectic group, we refer to \cite{paper3}. 


\section{Preliminaries on Clifford algebra}
\label{cliffordalgebra}


For a detailed description of the structure of Clifford algebras we refer to e.g.\ \cite{port}. Here we only recall the necessary basic notions. \\[-2mm]

The real Clifford algebra $\mathbb{R}_{0,m}$ is constructed over the vector space $\mathbb{R}^{0,m}$ endowed with a non--degenerate quadratic form of signature $(0,m)$, and generated by the orthonormal basis $(e_1,\ldots,e_m)$. The non--commutative Clifford or geometric multiplication in $\mathbb{R}_{0,m}$ is governed by the rules 
\begin{equation}\label{multirules}
e_{\alpha} e_{\beta} + e_{\beta} e_{\alpha} = -2 \delta_{\alpha \beta} \ \ , \ \ \alpha,\beta = 1,\ldots ,m
\end{equation}
As a basis for $\mathbb{R}_{0,m}$ one takes for any set $A=\lbrace j_1,\ldots,j_h \rbrace \subset \lbrace 1,\ldots,m \rbrace$ the element $e_A = e_{j_1} \ldots e_{j_h}$, with $1\leq j_1<j_2<\cdots < j_h \leq m$, together with $e_{\emptyset}=1$, the identity element. The dimension of $\mR_{0,m}$ is $2^m$. Any Clifford number $a$ in $\mathbb{R}_{0,m}$ may thus be written as $a = \sum_{A} e_A a_A$, $a_A \in \mathbb{R}$, or still as $a = \sum_{k=0}^m \lbrack a \rbrack_k$, where $\lbrack a \rbrack_k = \sum_{|A|=k} e_A a_A$ is the so--called $k$--vector part of $a$.
Real numbers thus correspond with the zero--vector part of the Clifford numbers. Euclidean space $\mathbb{R}^{0,m}$ is embedded in $\mathbb{R}_{0,m}$ by identifying $(X_1,\ldots,X_m)$ with the Clifford $1$--vector $\uX = \sum_{\alpha=1}^m e_{\alpha}\, X_{\alpha}$, for which it holds that $\uX^2 = - |\uX|^2 = - r^2$.\\[-2mm]

When allowing for complex constants, the generators $(e_1,\ldots, e_{m})$, still satisfying (\ref{multirules}), produce the complex Clifford algebra $\mathbb{C}_{m} = \mathbb{R}_{0,m} \oplus i\, \mathbb{R}_{0,m}$. Any complex Clifford number $\lambda \in \mathbb{C}_{m}$ may thus be written as $\lambda = a + i b$, $a,b \in \mathbb{R}_{0,m}$, leading to the definition of the Hermitian conjugation $\lambda^{\dagger} = (a +i b)^{\dagger} = \overline{a} - i \overline{b}$, where the bar notation stands for the Clifford conjugation in $\mathbb{R}_{0,m}$, i.e. the main anti--involution for which $\overline{e}_{\alpha} = -e_{\alpha}$, $\alpha=1, \ldots,m$. This Hermitian conjugation leads to a Hermitian inner product on $\mathbb{C}_{m}$ given by $(\lambda,\mu) = \lbrack \lambda^{\dagger} \mu \rbrack_0$ and its associated norm $|\lambda| = \sqrt{ \lbrack \lambda^{\dagger} \lambda \rbrack_0} = ( \sum_A |\lambda_A|^2 )^{1/2}$.\\[-2mm]

The algebra of real quaternions is denoted by $\mH$. For a quaternion
$$
q = q_0 + q_1 i + q_2 j + q_3 k = (q_0 + q_1 i )+ (q_2 + q_3 i)j = z + wj
$$
its conjugate is given by 
$$
\ol{q} = q_0 - q_1 i - q_2 j - q_3 k = (q_0 - q_1 i) - j(q_2  - q_3 i) = \ol{z} -j \ol{w} =  \ol{z} - w j
$$
such that
$$
q \ol{q} = \ol{q} q = |q|^2 = q_0^2 + q_1^2 + q_2^3 + q_3^2 = |z|^2 + |w|^2
$$
Identifying the quaternion units $i, j$ with the respective basis vectors $e_1, e_2$, the algebra $\mH$ is isomorphic with the Clifford algebra $\mR_{0,2}$. Moreover, it is also isomorphic with the even subalgebra $\mR_{0,3}^+$ of the Clifford algebra $\mR_{0,3}$, by identifying the quaternion units $i, j,k$ with the respective bivectors $e_2 e_3, e_3 e_1, e_1 e_2$.


\section{Euclidean and Hermitian Clifford Analysis: the basics}
\label{EHCA}


The central notion in standard Clifford analysis  is that of a {\em monogenic function}. This is a continuously differentiable function defined in an open region of Euclidean space $\mR^m$, taking its values in the Clifford algebra $\mR_{0,m}$, or subspaces thereof, and vanishing under the action of the Dirac operator $\dirac = \sum_{\alpha=1}^m e_{\alpha}\, \p_{X_{\alpha}}$, i.e.\ a vector valued first order differential operator, which can be seen as the Fourier or Fischer dual of the Clifford variable $\uX$. Monogenic functions thus are the higher dimensional counterparts of holomorphic functions in the complex plane. As moreover the Dirac operator factorizes the Laplacian: $\Delta_m = - \dirac^2$, the notion of monogenicity can be regarded as a refinement of the notion of harmonicity. It is important to note that the Dirac operator is invariant under the action of the $\mathrm{Spin}(m)$--group, which doubly covers the $\mathrm{SO}(m)$--group, whence this setting is usually referred to as Euclidean (or orthogonal) Clifford analysis.\\[-2mm]

Taking the dimension of the underlying Euclidean vector space $\mR^m$ to be even: $m=2n$, renaming the variables as:
$$
(X_1,\ldots,X_{2n}) = (x_1,y_1,x_2,y_2,\ldots,x_n,y_n)
$$
and considering the (almost) complex structure $\mI_{2n}$, i.e. the complex linear real $\mbox{SO}(2n)$ matrix 
$$
\mI_{2n} = \mathrm{diag} \begin{pmatrix} \phantom{-} 0 & 1 \\ -1 & 0 \end{pmatrix}$$
for which $\mI_{2n}^2 = - E_{2n}$, $E_{2n}$ being the identity matrix, we define the rotated vector variable
\begin{eqnarray*}
\uX_\mI = \mI_{2n}[\uX] & = & \mI_{2n} \left[ \sum_{k=1}^n (e_{2k-1} x_k + e_{2k} y_k ) \right]\\
& = & \sum_{k=1}^n \mI_{2n} [e_{2k-1}] x_k + \mI_{2n} [e_{2k}] y_k \ = \ \sum_{k=1}^n (-y_k e_{2k-1} + x_k e_{2k})
\end{eqnarray*}
and, correspondingly, the rotated Dirac operator
$$
\dirac_\mI = \mI_{2n} [ \dirac ] = \sum_{k=1}^n ( - \p_{y_k} e_{2k-1} + \p_{x_k} e_{2k})
$$
A differentiable function $F$ then is called {\em Hermitian monogenic} in some region $\Omega$ of $\mR^{2n}$, if and only if in that region $F$ is a solution of the system
\begin{equation}
\dirac F = 0 = \dirac_\mI F
\label{hmon}
\end{equation}

\noindent
Observe that this notion of Hermitian monogenicity does not involve the use of complex numbers, but could be developed as a real function theory instead. There is however an alternative approach to the concept of Hermitian monogenicity, making use of the projection operators $$\pi^{\mp} = \frac{1}{2} (\mp \, E_{2n}  + i \, \mI_{2n})$$ and thus involving a complexification.
In this approach we first define, in the complexification $\mC^{2n}$ of $\mR^{2n}$, the so--called Witt basis vectors
$$
\gf_k = \pi^- [e_{2k-1}] \quad \mathrm{and} \quad \gfd_k =  \pi^+ [e_{2k-1}] \quad (k=1,\ldots,n)
$$
for which also hold the relations
$$
i \, \gf_k = \pi^- [e_{2k}] \quad \mathrm{and} \quad -i \, \gfd_k =  \pi^+ [e_{2k}] \quad (k=1,\ldots,n)
$$
This enables the introduction of the vector variables $\uz$ and $\uzd$ given by
\begin{eqnarray*}
\uz & = & \pi^- [\uX] \ = \  \sum_{k=1}^n x_k \, \pi^- [e_{2k-1}] + \sum_{k=1}^n y_k \, \pi^- [ e_{2k} ] \\
&=& \sum_{k=1}^n \, x_k \gf_k + y_k (i \gf_k)  \ = \ \sum_{k=1}^n (x_k + i y_k) \, \gf_k = \sum_{k=1}^n \, z_k \, \gf_k \\
\uzd &=& \pi^+ [\uX] = \sum_{k=1}^n  x_k \gfd_k + y_k (-i \gfd_k)  = \sum_{k=1}^n \, \overline{z}_k \, \gfd_k
\end{eqnarray*}
and, correspondingly, of the Hermitian Dirac operators $\upz$ and $\upzd$ given by
\begin{eqnarray*}
2 \, \upzd &=& \pi^- [\dirac] \ = \ \sum_{k=1}^n \, \gf_k \p_{x_k} + i \gf_k \p_{y_k} = \sum_{k=1}^n \, \gf_k ( \p_{x_k} + i \p_{y_k} ) = 2 \, \sum_{k=1}^n \ \p_{\olz_k} \, \gf_k \\
2 \, \upz &=& \pi^+[\dirac] \ = \ \sum_{k=1}^n \, \gfd_k \p_{x_k} - i \gfd_k \p_{y_k} = \sum_{k=1}^n \, \gfd_k ( \p_{x_k} - i \p_{y_k} ) = 2 \, \sum_{k=1}^n \, \p_{z_k} \, \gfd_k
\end{eqnarray*}
As $2(\upz - \upzd) = \dirac$ and $2(\upz+\upzd) = i \, \dirac_\mI$, the system (\ref{hmon}) is easily seen to be equivalent with the system
\begin{equation}
\upz F = 0 = \upzd F
\label{hmoneq}
\end{equation}
Note that the $\cdot^\dagger$--notation corresponds to the Hermitian conjugation in the Clifford algebra $\mC_{2n}$, introduced above.\\[-2mm]

In order to study Hermitian monogenic functions it suffices to consider functions with values not in the whole Clifford algebra but in spinor space instead. Indeed, a Clifford algebra may be decomposed as a direct sum of isomorphic copies of a spinor space $\mS$, which, abstractly, may be defined as a minimal left ideal in the Clifford algebra. A spinor space is an irreducible Spin$(2n)$ group representation, and may be realized as follows. The Witt basis vectors satisfy the Grassmann identities
$$
\gf_j \gf_k + \gf_k \gf_j = 0 , \quad \gf_j^\dagger \gf_k^\dagger + \gf_k^\dagger \gf_j^\dagger = 0, \qquad j,k=1,\ldots,n
$$
which include their isotropy:
$$
\gf_j^2 = (\gf_j^\dagger)^2 = 0, \qquad j=1,\ldots, n
$$
as well as the duality identities
$$
\gf_j \gf_k^\dagger + \gf_k^\dagger \gf_j = \delta_{jk}, \qquad j,k=1,\ldots,n
$$
The Witt basis vectors $(\gf_1,\ldots,\gf_{n})$ and $(\gf_1^\dagger,\ldots,\gf_{n}^\dagger)$ respectively span isotropic subspaces $W$ and $W^\dagger$ of $\mC^{2n}$, such that 
\begin{equation}
\label{isotropic}
\mC^{2n} = W \oplus W^\dagger
\end{equation}
those subspaces being eigenspaces of the complex structure $\mI_{2n}$ with respective eigenvalues $-i$ and $i$. They also generate two Grassmann algebras, respectively denoted by $\mC \bigwedge_{n}$ and $\mC \bigwedge_{n}^\dagger$. \\[-2mm]

With the self-adjoint idempotents 
$$
I_j = \gf_j \gf_j^\dagger = \frac{1}{2} ( 1- i e_{2j-1} e_{2j} ), \qquad j=1,\ldots,n
$$
we compose the primitive self--adjoint idempotent $I = I_1 I_2 \cdots I_{n}$ leading to the realization of the spinor space $\mS$ as $\mS = \mC_{2n} I$. Since $\gf_j I =0$, $j=1,\ldots,n$, we also have $\mS \simeq \mC \bigwedge_{n}^\dagger I$.\\[-2mm]

\noindent
When decomposing the Grassmann algebra $\mC \bigwedge_{n}^\dagger$ into its so--called homogeneous parts
$$
\mC \bigwedge_{n}\nolimits^\dagger = \bigoplus_{r=0}^{n} \biggl ( \mC \bigwedge_{n}\nolimits^\dagger \biggr )^{(r)}
$$
where the spaces $\left (\mC \bigwedge_{n}^\dagger \right )^{(r)}$ are spanned by all possible products of $r$ Witt basis vectors out of $(\gf_1^\dagger,\ldots,\gf_{n}^\dagger)$, the spinor space $\mS$ accordingly decomposes into
\begin{equation}
\mS = \bigoplus_{r=0}^{n} \mS^r, \qquad \mbox{with \ } \mS^r \simeq \left ( \mC \bigwedge_{n}\nolimits^\dagger \right)^{(r)} I
\label{decompspin}
\end{equation}
These homogeneous parts $\mS^r$, $r=0,\ldots,n$, provide models for fundamental $\mbox{U}(n)$--representations (see \cite{howe}) and thus also for fundamental $\gsl(n,\mC)$-representations (see \cite{partI}, \cite{Ee}). \\[-2mm]

Accordingly we can decompose a spinor valued function $F: \mC^{n} \longrightarrow \mS$ into its components 
$$
F = \sum_{r=0}^n F^r, \quad F^r: \mC^{n} \longrightarrow \mS^r, \quad r=0,\ldots,n
$$
Pay attention to the fact that the monogenicity of $F$ does {\em not} imply the monogenicity of the components $F^r$, $r=0,\ldots,n$. However, the Hermitian monogenicity of $F$ does imply the Hermitian monogenicity of these components, and vice versa. This is due to the nature of the action of the Witt basis vectors as (left) multiplication operators, implying that
$$
\upz F^r: \mC^{n} \longrightarrow \mS^{r+1} \quad \mbox{and} \quad \upzd F^r: \mC^{n} \longrightarrow \mS^{r-1} 
$$
Moreover, for each of the components $F^r$ the notions of monogenicity and Hermitian monogenicity coincide, since
$$
\dirac F^r = 0 \Longleftrightarrow (\upz - \upzd) F^r = 0 \Longleftrightarrow \left \{ \begin{array}{lcl}
\upz F^r & = & 0 \\[1mm] \upzd F^r & = & 0 \end{array} \right .
$$
In conclusion we have the following result.

\begin{proposition}
The function $F = \sum_{r=0}^n F^r$, $F^r : \mC^{n} \longrightarrow \mS^r$, $r=0,\ldots,n$, is Hermitian monogenic in some region $\Omega$ of $\mC^{n}$ if and only if each of its components $F^r$, $r=0,\ldots,n$, is monogenic in that region.
\end{proposition}

The Hermitian Dirac operators $\upz$ and $\upzd$ are invariant under the action of the group $\mbox{SO}_{\mI}(2n)$, i.e.\ the subgroup of $\mbox{SO}(2n)$ consisting of those matrices which commute with the complex structure $\mI_{2n}$. This subgroup $\mbox{SO}_\mI(2n)$ inherits a twofold covering by the subgroup $\mbox{Spin}_{\mI}(2n)$ of $\mbox{Spin}(2n)$, consisting of those elements of $\mbox{Spin}(2n)$ that are commuting with 
\begin{equation}
s_{\mI} = s_1 \ldots s_{n},  \quad \mbox{where \ } s_j = \frac{\sqrt{2}}{2} ( 1 + e_{2j-1} e_{2j} ),  \; j=1,\ldots,n
\label{sI}
\end{equation}
This element $s_{\mI}$ itself obviously belongs to $\mbox{Spin}_{\mI}(2n)$ and corresponds, under the double covering, to the complex structure $\mI_{2n} $. As $\mbox{SO}_{\mI}(2n)$ is isomorphic with the unitary group $\mbox{U}(n)$ and, up to this isomorphism, $\mbox{Spin}_{\mI}(2n)$ thus provides a double cover of $\mbox{U}(n)$, we may say that $\mbox{U}(n)$ is the fundamental group underlying the function theory of Hermitian monogenic functions. 


\section{Quaternionic Clifford Analysis}
\label{QCA}


A refinement of Hermitian Clifford analysis is obtained by considering the hypercomplex structure $\mQ = (\mI_{4p},\mJ_{4p},\mK_{4p})$ on $\mR^{4p} \simeq \mC^{2p} \simeq \mH^p$, the dimension $m=2n=4p$ now being assumed to be a multiple of four. This hypercomplex structure arises by introducing, next to the first complex structure $\mI_{4p}$, a second one, called $\mJ_{4p}$, given by
$$
\mJ_{4p} = \mathrm{diag} \, \left ( \begin{array}{cccc} &  & 1 & \\ & & & -1 \\  -1 & & &  \\ & 1 & & \end{array} \right )
$$
Clearly $\mJ_{4p}$ belongs to $\mbox{SO}(4p)$, with $\mJ_{4p}^2 = -E_{4p}$, and it anti--commutes with $\mI_{4p}$. 
Then a third complex structure quite naturally arises, namely the $\mbox{SO}(4p)$--matrix 
$$
\mK_{4p} = \mI_{4p} \, \mJ_{4p} = - \mJ_{4p} \, \mI_{4p}
$$ 
for which $\mK_{4p}^2 = -E_{4p}$ and which anti--commutes with both $\mI_{4p}$ and $\mJ_{4p}$. It turns out that 
$$
\mK_{4p} = \mbox{diag} \, \left ( \begin{array}{cccc} &  &  & -1 \\ & & -1 &  \\   & 1 & &  \\  1 & & & \end{array} \right )
$$
The $\mbox{SO}(4p)$--matrices which commute with the hypercomplex structure $\mQ$ on $\mR^{4p}$ form a subgroup of $\mbox{SO}_{\mI}(4p)$, denoted by $\mbox{SO}_{\mQ}(4p)$, which is isomorphic with the symplectic group $\mbox{Sp}( p)$.
Recall that the symplectic group $\mbox{Sp}(p)$ is the real Lie group of quaternion $p \times p$ matrices preserving the symplectic inner product
$$
\langle \xi , \eta \rangle_{\mH} = \xi_1 \overline{\eta}_1 + \xi_2 \overline{\eta}_2 + \cdots +  \xi_p \overline{\eta}_p \quad \xi, \eta \in \mH^p
$$
or, equivalently,
$$
\mbox{Sp}(p) = \left \{ A \in \mbox{GL}_p(\mH) : AA^\ast = E_p \right \}
$$  
Quite naturally, the subgroup $\mbox{SO}_{\mQ}(4p)$ of $\mbox{SO}(4p)$ is doubly covered by $\mbox{Spin}_{\mQ}(4p)$, the subgroup of $\mbox{Spin}(4p)$ consisting of the $\mbox{Spin}(4p)$--elements which are commuting with both $s_{\mI}$ and $s_{\mJ}$, where now $s_{\mJ}$ is the $\mbox{Spin}(4p)$--element corresponding to $\mJ_p$. Recall, see (\ref{sI}), that $s_\mI$, corresponding to the complex structure $\mI_{4p}$, is given by $s_{\mI} = s_1 \cdots s_{2p}$, where $s_j = \frac{\sqrt{2}}{2} \bigl ( 1 + e_{2j-1} e_{2j} \bigr )$, $j=1,\ldots, 2p$. Similarly, for $s_{\mJ}$ we find
\begin{equation}
s_{\mJ} = \widetilde{s_1} \cdots \widetilde{s_p}, \quad \widetilde{s_j} = \frac{1}{2} \bigl ( 1 + e_{4j-3} e_{4j-1} \bigr ) \bigl (1 - e_{4j-2} e_{4j} \bigr ), \; j=1,\ldots, p
\label{sJ}
\end{equation}
For the corresponding picture at the level of the Lie algebras we refer to \cite{paper1}.\\[-2mm]

The introduction of a hypercomplex structure leads to a function theory in the framework of so--called {\em quaternionic Clifford analysis}, where the fundamental invariance will be that of the symplectic group $\mbox{Sp}(p)$. The most genuine way to introduce the new concept of quaternionic monogenicity is to directly generalize the system (\ref{hmon}) expressing Hermitian monogenicity, now making use of the hypercomplex structure on $\mR^{4p}$ and the additional rotated Dirac operators $\dirac_\mJ = \mJ_{4p}[\dirac]$ and $\dirac_\mK = \mK_{4p}[\dirac]$; whence the following definition.

\begin{definition}
\label{defqmon}
A differentiable function $F: \mR^{4p} \longrightarrow \mS$ is called quaternionic monogenic (q--monogenic for short) in some region $\Omega$ of $\mR^{4p}$, if and only if in that region $F$ is a solution of the system
\begin{equation}
\dirac F = 0, \quad \dirac_\mI F = 0, \quad \dirac_\mJ F = 0, \quad \dirac_\mK F = 0
\label{qmon}
\end{equation}
\end{definition}
\noindent
Observe that, in a similar way as it was possible to introduce the notion of Hermitian monogenicity without having to resort to complex numbers, the above Definition \ref{defqmon} expresses the notion of q--monogenicity without  involving quaternions.

\begin{remark}
The notion of a hypercomplex structure stems from differential geometry where each tangent bundle of an even dimensional manifold admits the action of the algebra of quaternions, the quaternion units defining three almost complex structures. Accordingly we should in fact call the functions satisfying system (\ref{qmon}) ``hypercomplex monogenic'', however this term has already been used  in other contexts in higher dimensional function theory, and moreover, the notion of quaternionic monogenicity was already introduced in previous papers, whence we stick to the latter terminology.

\end{remark}

There is a natural alternative characterization of q--monogenicity possible in terms of the Hermitian Dirac operators, yet still not involving quaternions. We recall these Hermitian Dirac operators in the current dimension:
\begin{eqnarray*}
\upz & = & \sum_{k=1}^{2p} \p_{z_k} \gfd_k \ = \ \sum_{j=1}^p ( \p_{z_{2j-1}} \gfd_{2j-1} + \p_{z_{2j}} \gfd_{2j} ) \\
\upzd & = & \sum_{k=1}^{2p} \p_{\overline{z}_k} \gf_k \ = \ \sum_{j=1}^p ( \p_{\overline{z}_{2j-1}} \gf_{2j-1} + \p_{\overline{z}_{2j}} \gf_{2j} )
\end{eqnarray*}
and compute their respective images under the action of the complex structure $\mJ_{4p}$:
\begin{eqnarray*}
\upzJ \ = \ \mJ_{4p} [\upz] & = &  \sum_{j=1}^p ( \p_{z_{2j}} \gf_{2j-1} - \p_{z_{2j-1}} \gf_{2j} ) \\
\upzJd \ = \ \mJ_{4p} [\upzd] & = & \sum_{j=1}^p ( \p_{\overline{z}_{2j}} \gfd_{2j-1} - \p_{\overline{z}_{2j-1}} \gfd_{2j} )
\end{eqnarray*}
Similarly we can introduce the auxiliary variables
\begin{eqnarray*}
\uzJ \ = \ \mJ_{4p} [\uz] & = & \sum_{j=1}^p ( z_{2j} \, \gfd_{2j-1} - z_{2j-1} \, \gfd_{2j} ) \\
\uzJd \ = \ \mJ_{4p} [\uzd] & = & \sum_{j=1}^p ( \overline{z}_{2j} \, \gf_{2j-1} - \overline{z}_{2j-1} \, \gf_{2j} )
\end{eqnarray*}
Here, use has been made of the formulae
$$
\mJ_{4p} [ \gf_{2j-1} ] = - \gfd_{2j}, \quad \mJ_{4p} [\gf_{2j}] = \gfd_{2j-1}, \quad \mJ_{4p}[\gfd_{2j-1}] = - \gf_{2j}, \quad \mbox{and} \quad \mJ_{4p}[\gfd_{2j}] = \gf_{2j-1}
$$
Now the original Dirac operator $\dirac$ and its rotated versions $\dirac_\mI$, $\dirac_\mJ$ and $\dirac_\mK$ may be expressed in terms of the Hermitian Dirac operators $( \upz,\upzd )$ and their $\mJ$--rotated versions $( \upzJ,\upzJd )$ as follows:
\begin{eqnarray*}
\upz &=& \pi^+[\dirac] = \frac{1}{4} ( {\bf 1} + i \, \mI_{4p} ) [\dirac], \quad \upzd = \pi^-[\dirac] = -\frac{1}{4} ( {\bf 1} - i \, \mI_{4p} ) [\dirac] \\
\upzJ &=& \mJ_{4p} [\pi^+[\dirac]] = \frac{1}{4} ( \mJ_{4p} + i \, \mK_{4p}) [\dirac], \quad \upzJd =  \mJ_{4p} [\pi^-[\dirac]] = -\frac{1}{4} ( \mJ_{4p} - i \, \mK_{4p}) [\dirac]
\end{eqnarray*}
whence conversely
\begin{eqnarray*}
\dirac & = & 2 ( \upz - \upzd) \\
i \, \dirac_\mI & = & 2 ( \upz + \upzd ) \\
\dirac_\mJ &=& 2 ( \upzJ - \upzJd ) \\
i \, \dirac_\mK &=& 2 ( \upzJ + \upzJd) 
\end{eqnarray*}
Clearly, this leads to an alternative characterization of q--monogenicity as expressed in the following proposition.
\begin{proposition}
A differentiable function $F: \mR^{4p} \simeq \mC^{2p} \longrightarrow \mS$ is q--monogenic in some region $\Omega \subset \mR^{4p}$ if and only if $F$ is in $\Omega$ a simultaneous null solution of the operators $\upz$, $\upzd$, $\upzJ$ and $\upzJd$.
\end{proposition}

We want to emphasize that, in Section \ref{gradients}, the new operators $\upzJ$ and $\upzJd$, used, together with the known Hermitian Dirac operators $\upz$ and $\upzd$, to characterize the concept of q--monogenicity, will arise in a natural way as generalized gradients in the sense of Stein \& Weiss (see \cite{stein}).\\[-2mm]

As the identification of an underlying symmetry group is necessary for the further development of the function theory, the following result is crucial.
\begin{proposition}
The operators $\upz$, $\upzd$, $\upzJ$ and $\upzJd$ are invariant under the action of the symplectic group $\mbox{\em Sp}(p)$.
\end{proposition}

\pf
The action of a $\mbox{Spin}(4p)$--element $s$ on a spinor valued function $F$ is the so-called $L$--action given by $L(s) [ F(\uX)] = s F ( s^{-1} \uX s)$. The Dirac operator $\dirac$ is invariant under $\mbox{Spin}(4p)$, i.e.
$$
[L(s),\dirac] = 0, \qquad \mbox{for all } s \in \mbox{Spin}(4p)
$$
which can be explained by
$$
L(s) \dirac_{\uX} F(\uX) = s \dirac_{s^{-1} \uX s} F (s^{-1} \uX s) = s ( s^{-1} \dirac_{\uX} s ) F ( s^{-1} \uX s)  = \dirac_{\uX} L(s) F(\uX)
$$
Recall that $\mbox{Sp}(p)$ is isomorphic to the subgroup Spin$_{\mQ}(4p)$ of $\mbox{Spin}(4p)$, whence the Dirac operator $\dirac$ is, quite trivially, also invariant under the action of $\mbox{Sp}(p)$. The invariance of the operators $\upz$, $\upzd$, $\upzJ$ and $\upzJd$ now follows from the fact that their respective definitions only involve projection operators which are commuting with the $\mbox{Sp}(p)$--elements.
\qed

We will now comment on the behaviour of the notion of q--monogenicity with respect to the decomposition of a function $F = \sum_{r=0}^n F^r$, $F^r : \mR^{4p} \longrightarrow \mS^r$, $r=0,\ldots,n$ into its homogeneous spinor components. \\[-2mm]

If $F$ is q--monogenic, then so are its components $F^r$, and vice versa. If $F^r$ is q--monogenic, then, quite naturally, $F^r$ is Hermitian monogenic with respect to the Hermitian Dirac operators $\upz$ and $\upzd$ --- let us call this ${\mI}$--Hermitian monogenicity --- but also Hermitian monogenic with respect to the rotated Hermitian Dirac operators $\upzJ$ and $\upzJd$ --- let us call this $\mJ$--Hermitian monogenicity. As was already pointed out, for the function $F^r$, ${\mI}$--Hermitian monogenicity is equivalent with $\dirac$--monogenicity; in the same order of ideas, the $\mJ$--Hermitian monogenicity of $F^r$ is equivalent with its $\dirac_\mJ$--monogenicity. Summarizing, we have the following result.

\begin{proposition}
\label{homog-quaternionic}
For the function $F = \sum_{r=0}^n F^r$, $F^r : \mR^{4p} \longrightarrow \mS^r$, $r=0,\ldots,n$, defined in some region $\Omega$ of $\mR^{4p}$, the following statements are equivalent:
\begin{itemize}
\item[(i)] $F$ is q--monogenic;\\[-7mm]
\item[(ii)] each of the components $F^r$ is q--monogenic;\\[-7mm]
\item[(iii)] each of the components $F^r$ is simultaneously $\dirac$--monogenic and $\dirac_\mJ$--monogenic.
\end{itemize}
\end{proposition}


\section{A further decomposition of spinor space}
\label{spinordecomposition}


Spinor space $\mS$, which was already decomposed into U$(2p)$--irreducible homogeneous parts $\mS^r$, can further be decomposed into Sp$( p)$--irreducibles, which we will call {\em symplectic cells}. Let us briefly sketch this decomposition, referring to \cite{paper1} for a detailed description.\\[-2mm]

First we introduce the Sp$( p)$--invariant left multiplication operators
\begin{eqnarray*}
P &=& \gf_2 \gf_1 + \gf_4 \gf_3 + \ldots + \gf_{2p} \gf_{2p-1} \\
Q &=& \gfd_1 \gfd_2 + \gfd_3 \gfd_4 + \ldots + \gfd_{2p-1} \gfd_{2p} \ = \ P^{\dagger}
\end{eqnarray*}
for which $P: \mS^r \rightarrow \mS^{r-2}$ and $Q: \mS^r \rightarrow \mS^{r+2}$. 
Next we define the symplectic cells of spinor space, which are given by, for $r=0,\ldots,p$, the subspaces
$$
\mS_r^r = \mbox{Ker} \, P |_{\mS^r}, \qquad \mS_r^{2p-r} = \mbox{Ker} \, Q |_{\mS^{2p-r}}
$$
and for $k=0,\ldots,p-r$, the subspaces
$$
\mS_r^{r+2k} = Q^k \, \mS^r_r, \qquad \mS_r^{2p-r-2k} = P^k \, \mS^{2p-r}_r
$$

\begin{lemma}
\label{lemmaS4}
One has
$$
\mbox{\em Ker} \, P|_{\mS} = \bigoplus_{r=0}^p \mS_r^r \qquad \mbox{and} \qquad \mbox{\em Ker} \, Q|_{\mS} = \bigoplus_{r=0}^p \mS_r^{2p-r} 
$$ 
and, for $k= 0,1,\ldots,p-r-1$,
\begin{itemize}
\item $Q$ is an isomorphism $\mS_r^{r+2k} \longrightarrow \mS_r^{r+2k+2}$ with inverse $Q^{-1} = \displaystyle\frac{1}{\alpha_r^k} \, P$;\\[-4mm]
\item  $P$ is an isomorphism $\mS_r^{2p-r-2k} \longrightarrow \mS_r^{2p-r-2k-2}$ with inverse $P^{-1} = \displaystyle\frac{1}{\alpha_r^{p-r-k-1}} \, Q$
\end{itemize}
where the coefficients are given by
$$
\alpha_r^k = (k+1)(p-r-k) = \alpha_r^{p-r-k-1}
$$
which implies that the composition of the multiplicative operators $P$ and $Q$ is constant on each symplectic cell; more specifically one has
\begin{itemize}
\item $P \, Q = \alpha_r^k$ on \, $\mS_r^{r+2k}$ and on \, $\mS_r^{2p-r-2k-2}$
\item $Q \, P = \alpha_r^k$ on \,  $\mS_r^{r+2k+2}$ and on \, $\mS_r^{2p-r-2k}$
\end{itemize}
\end{lemma}

\begin{proposition}
\label{propositionS1}
One has, for all $r=0,\ldots,p$,
$$
\mS^r = \bigoplus_{j=0}^{\lfloor \frac{r}{2} \rfloor} \mS_{r-2j}^r \qquad \mbox{and} \qquad
\mS^{2p-r} = \bigoplus_{j=0}^{\lfloor \frac{r}{2} \rfloor} \mS_{r-2j}^{2p-r}
$$
and each of the symplectic cells $\mS_s^r$ in the above decompositions is an irreducible $\mbox{Sp}(p)$--representation. 
\end{proposition}

\noindent
The above decomposition of the homogeneous spinor subspaces into symplectic cells can be depicted by the following triangular scheme.
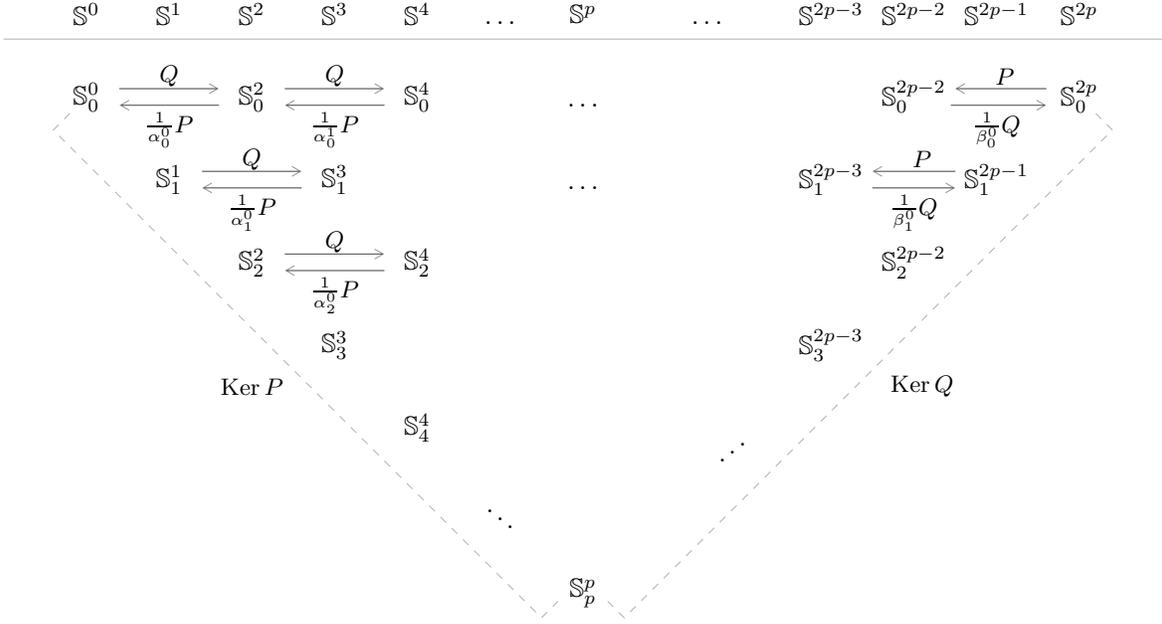
\begin{center}
\begin{tikzpicture}[scale=1.1]
\node at (0,0) {$\mS^0$};
\node at (1,0) {$\mS^1$};
\node at (2,0) {$\mS^2$};
\node at (3,0) {$\mS^3$};
\node at (4,0) {$\mS^4$};
\node at (5,-0.1) {$\ldots$};
\node at (6,0) {$\mS^p$};
\node at (7.5,-0.1) {$\ldots$};
\node at (9,0) {$\mS^{2p-3}$};
\node at (10,0) {$\mS^{2p-2}$};
\node at (11,0)  {$\mS^{2p-1}$};
\node at (12,0) {$\mS^{2p}$};
\draw[-, draw opacity=0.3] (-1,-0.3) -- (13,-0.3);
\node at (0,-1) {$\mS^0_0$};
\node at (2,-1) {$\mS^2_0$};
\node at (4,-1) {$\mS^4_0$};
\node at (6,-1.1) {$\ldots$};
\node at (6,-2.1) {$\ldots$};
\node at (10,-1) {$\mS^{2p-2}_0$};
\node at (12,-1) {$\mS^{2p}_0$};
\node at (1,-2) {$\mS^1_1$};
\node at (3,-2) {$\mS^3_1$};
\node at (11,-2) {$\mS^{2p-1}_1$};
\node at (9,-2) {$\mS^{2p-3}_1$};
\node at (2,-3) {$\mS_2^2$};
\node at (4,-3) {$\mS_2^4$};
\node at (10,-3) {$\mS_2^{2p-2}$};
\node at (9,-4) {$\mS_3^{2p-3}$};
\node at (3,-4) {$\mS_3^3$};
\node at (4,-5) {$\mS_4^4$};
\node at (5,-6) {$\ddots$};
\node at (6,-7) {$\mS_p^p$};
\node at (7.8,-5.2) {$\iddots$};
\draw[->,>=angle 60,draw opacity=0.4] (0.4,-0.9) -- (1.6,-0.9);
\node at (1,-0.75) {\small $Q$};
\draw[<-,>=angle 60,draw opacity=0.4] (0.4,-1.1) -- (1.6,-1.1);
\draw[<-,>=angle 60,draw opacity=0.4] (10.5,-0.9) -- (11.6,-0.9);
\draw[->,>=angle 60,draw opacity=0.4] (10.45,-1.1) -- (11.6,-1.1);
\node at (11.1,-0.75) {\small $P$};
\node at (3,-0.75) {\small $Q$};
\draw[->,>=angle 60,draw opacity=0.4] (2.4,-0.9) -- (3.6,-0.9);
\draw[<-,>=angle 60,draw opacity=0.4] (2.4,-1.1) -- (3.6,-1.1);
\node at (2,-1.75) {\small $Q$};
\draw[->,>=angle 60,draw opacity=0.4] (1.4,-1.9) -- (2.6,-1.9);
\draw[<-,>=angle 60,draw opacity=0.4] (1.4,-2.1) -- (2.6,-2.1);
\draw[<-,>=angle 60,draw opacity=0.4] (9.5,-1.9) -- (10.5,-1.9);
\draw[->,>=angle 60,draw opacity=0.4] (9.5,-2.1) -- (10.5,-2.1);
\node at (10.1,-1.75) {\small $P$};
\node at (3,-2.75) {\small $Q$};
\draw[->,>=angle 60,draw opacity=0.4] (2.4,-2.9) -- (3.6,-2.9);
\draw[<-,>=angle 60,draw opacity=0.4] (2.4,-3.1) -- (3.6,-3.1);
\draw[draw opacity = 0.3,dashed] (-0.2,-1.2) -- (-0.4,-1.4) -- (5.5,-7.3) -- (5.7,-7.1) ;
\node at (2,-4.5) {\small $\mbox{Ker} \, P$};
\draw[draw opacity = 0.3,dashed] (6.3,-7.1) -- (6.5,-7.3) -- (12.4,-1.4) -- (12.2,-1.2);
\node at (10.1,-4.5) {\small $\mbox{Ker} \, Q$};
\node at (1,-1.4) {\small $\frac{1}{\alpha_0^0} P$};
\node at (3,-1.4) {\small $\frac{1}{\alpha_0^1} P$};
\node at (2,-2.4) {\small $\frac{1}{\alpha_1^0} P$};
\node at (3,-3.4) {\small $\frac{1}{\alpha_2^0} P$};
\node at (11,-1.4) {\small $\frac{1}{\beta_0^0} Q$};
\node at (10,-2.4) {\small $\frac{1}{\beta_1^0} Q$};
\end{tikzpicture}
\end{center}

Now let us decompose a function $F: \mC^{2p} \longrightarrow \mS$ according to these symplectic cells of spinor space:
$$
F = \sum_{r=0}^{2p} F^r = \sum_{r=0}^{2p} \sum_s F_s^r, \qquad F^r: \mC^{2p} \longrightarrow \mS^r, \quad F_s^r : \mC^{2p} \longrightarrow \mS_s^r
$$
and investigate the possible inheritance of the various concepts of monogenicity by the distinguished components.

\begin{enumerate}
\item[(i)] If $F$ is monogenic, then nothing can be said about the monogenicity of the components $F^r_s$. This is because, a fortiori, the action of the Dirac operator on the function values
$$
\dirac = 2 (\upzd - \upz) : \mS^r \longrightarrow \mS^{r+1} \oplus \mS^{r-1}
$$
mixes up the homogeneous parts of spinor space.

\item[(ii)] If $F$ is Hermitian monogenic, then each symplectic component $F_s^r$ is not necessarily Hermitian monogenic; this because the action of each of the Hermitian Dirac operators on the function values is mixing up the symplectic cells of a homogeneous spinor subspace, since indeed
$$
\upz :  \mS^r_s \longrightarrow \mS^{r+1}_{s-1} \oplus \mS^{r+1}_{s+1}
\quad {\rm and} \quad
\upzd :  \mS^r_s \longrightarrow \mS^{r-1}_{s-1} \oplus \mS^{r-1}_{s+1}
$$

\item[(iii)] If $F$ is q--monogenic, then, remarkably and in contrast with Hermitian monogenicity, the symplectic components $F_s^r$ will  be q--monogenic, as we will prove in the next proposition.
\end{enumerate}

\begin{lemma}
\label{PQF}
If the function $F : \mC^{2p} \longrightarrow \mS$ is q--monogenic, then so are the functions $PF$ and $QF$.
\end{lemma}

\pf
This result directly follows from the commutation relations of the operators $P$ and $Q$ with the quaternionic Hermitian Dirac operators:
\begin{align}
\label{commutPQ}
[P , \upz] &=  - \upzJ  &[Q , \upz] &= \phantom{-} 0\nonumber\\
[P , \upzd] &= \phantom{-} 0   &[Q , \upzd] &=  \phantom{-} \upzJd\nonumber\\
[P , \upzJ] &= \phantom{-} 0  &[Q , \upzJ] &= - \upz\nonumber\\
[P , \upzJd] &=  \phantom{-}  \upzd  & [Q , \upzJd] &=  \phantom{-} 0
\end{align}
\eop

\begin{proposition}
\label{comp}
A~function $F : \mC^{2p} \longrightarrow \mS$ is q--monogenic if and only if its symplectic components $F^r_s$ are.
\end{proposition}

\pf We only need to prove that, for a given q--monogenic function $F : \mC^{2p} \longrightarrow \mS$, its symplectic components $F^r_s$ are q--monogenic as well. The opposite implication is obvious. Now, by Proposition \ref{homog-quaternionic}, we know that each component $F^r$ is q--monogenic. Putting $r=2t$, we have by Lemma \ref{PQF} that $F_0^{2t}\doteq Q^tP^tF^{2t}$ is q--monogenic. Here $F\doteq G$ means that $F$ and $G$ are equal up to a~non-zero multiple. Further, we have that $F_2^{2t}\doteq Q^{t-1}P^{t-1}(F^{2t}-F_0^{2t})$ is q--monogenic and so on.
\eop\\

As already pointed out above, when restricting the values of the functions considered to a homogeneous subspace $\mS^r$ of spinor space, Hermitian monogenicity can be expressed using only one operator, namely the Dirac operator. For functions taking their values in a symplectic cell $\mS^r_s$ a similar, quite remarkable, result is valid.
\begin{proposition}
\label{equivalence-q-monog}
For functions $F_s^r: \mC^{2p} \longrightarrow \mS_s^r$, the equations $\upz F_s^r=0$ and $\upzJ F_s^r=0$ are equivalent as long as $0 \leq s < r$, and the same holds for the equations $\upzd F_s^r=0$ and $\upzJd F_s^r=0$. Both equivalences also remain valid for $(r,s)=(p,p)$. If $0 \leq r < p$ and $s=r$ then $\upz F_r^r =0$ implies $\upzJ F_r^r=0$ and $\upzJd F_r^r=0$ implies $\upzd F_r^r=0$, but not conversely. If $p < r \leq 2p$ and $s=2p-r$ then $\upzJ F_{2p-r}^r = 0$ implies $\upz F_{2p-r}^r=0$ and $\upzd F_{2p-r}^r =0$ implies $\upzJd F_{2p-r}^r=0$, but not conversely.
\end{proposition}

\pf
We start with the exceptional case $F_p^p: \mC^{2p} \longrightarrow \mS^p_p$, the value space $\mS_p^p$ then being in the kernel of both operators $P$ and $Q$. Assume that $\upz F_p^p=0$, then $P \upz F_p^p=0$ or $( \upz P - \upzJ) F_p^p=0$, whence $\upzJ F_p^p=0$. Conversely, if $\upzJ F_p^p=0$, then $Q \upzJ F_p^p = 0$ or $(\upzJ Q - \upz) F_p^p =0$, whence $\upz F_p^p=0$. A similar reasoning shows that $\upzd F_p^p=0$ and $\upzJd F_p^p=0$ are equivalent systems of equations on $\mS^p_p$.\\
In general, for functions $F_s^r: \mC^{2p} \longrightarrow \mS_s^r$, we consecutively have
$$
\upz F_s^r = 0 \Rightarrow P \upz F_s^r = 0 \Rightarrow (\upz P - \upzJ ) F_s^r=0 \Rightarrow \upz P F_s^r =  \upzJ F_s^r
$$ 
Applying the operator $Q$ leads to $Q \upz P F_s^r =  Q \upzJ F_s^r$ or $Q \upzJ F_s^r =  \upz QP F_s^r$. However on $\mS_s^r$ the product $QP$ acts as a constant, depending upon $r$ and $s$, whence it follows that $Q \upzJ F_s^r = 0$. As long as $\upzJ F_s^r$ does not take values in $\mbox{Ker} \, Q = \bigoplus_{s=0}^p \mS_s^{2p-s}$ it follows that $\upzJ F_s^r=0$. In a similar way it is shown that $\upzJd F_s^r = 0$ implies $\upzd F_s^r = 0$, as long as $\upzd F_s^r$ does not take values in $\mbox{Ker} \, Q$. A reasoning along similar lines shows that $\upzJ F_s^r=0$ implies $\upz F_s^r=0$ and $\upzd F_s^r=0$ implies $\upzJd F_s^r=0$, as long as neither $\upz F_s^r$ nor $\upzJd F_s^r$ take values in $\mbox{Ker} \, P = \bigoplus_{s=0}^p \mS^s_s$.
\qed

\begin{corollary}
If the function $F^r_s$ takes its values in $\mS^r_s$ with $r=0,1,2,\ldots,2p$ and $0 \leq s <r$, the following statements are equivalent:
\begin{itemize}
\item[(i)] $F^r_s$ is $\dirac$--monogenic;\\[-7mm]
\item[(ii)] $F^r_s$ is $\dirac_\mJ$--monogenic;\\[-7mm]
\item[(iii)] $F^r_s$ is q--monogenic.
\end{itemize}
The same equivalence holds in the case where $(r,s)=(p,p)$. However, for $s=r \neq p$ this equivalence reduces to
\begin{itemize}
\item[(i')] $F^r_s$ is $\dirac$ and $\dirac_\mJ$--monogenic;\\[-7mm]
\item[(ii')] $F^r_s$ is q--monogenic.
\end{itemize}
\end{corollary}

\begin{remark}
\label{partial result}
In \cite{paper1} we introduced a decomposition of the multiplicative operators
$\gf_j$ and $\gfd_j$, $j=1,\ldots,p$:
\begin{eqnarray*}
\gf_j \Bigr\rvert_{\mS_s^r} &=& (\gf_j)^r_{s-} +  (\gf_j)^r_{s+} \quad , \quad (\gf_j)^r_{s\mp} : \mS_s^r \longrightarrow \mS_{s\mp1}^{r-1} \\
\gfd_j \Bigr\rvert_{\mS_s^r} &=& (\gfd_j)^r_{s-} +  (\gfd_j)^r_{s+} \quad , \quad (\gfd_j)^r_{s\mp} : \mS_s^r \longrightarrow \mS_{s\mp1}^{r+1}
\end{eqnarray*}
Now, we can decompose the four operators expressing q--monogenicity accordingly. If the function $F^r_s$ takes its values in $\mS^r_s$, then $(\upz)^r_{s\mp} \, F^r_s$ will take its values in $\mS_{s\mp1}^{r+1}$, and similarly for the other three operators. From the above Proposition \ref{equivalence-q-monog} we know that for $0 \leq r < p$, $\upzJ F_r^r=0$ does not imply $\upz F_r^r=0$, nor does $\upzd F^r_r = 0$ imply $\upzJd F^r_r = 0$. However we have a partial result:  it can be proven that $\upzJ F_r^r=0$ implies $(\upz)^r_{r-} F_r^r=0$ and that $\upzd F^r_r = 0$ implies $(\upzJd)^r_{s-} F^r_r = 0$. Similar results hold for the function $F^r_{2p-r}$ with $p < r \leq 2p$: $\upz F^r_{2p-r} = 0$ implies $(\upzJ)^r_{(2p-r)-} F^r_{2p-r} =0$ and $\upzJd F^r_{2p-r} = 0$ implies $(\upzd)^r_{(2p-r)-} F^r_{2p-r} =0$
\end{remark}

In the same order of ideas as in Remark \ref{partial result}, we can express the equivalences described in Proposition \ref{equivalence-q-monog} by means of the operators obtained through the decomposition of the multiplicative action of the Witt basis vectors. The equivalence of the operators $\upz$ and $\upzJ$ incorporates both the equivalence of the operators $(\upz)_-$ and $(\upzJ)_-$ and the equivalence of the operators $(\upz)_+$ and $(\upzJ)_+$. Similarly, the equivalence of  the operators $\upzd$ and $\upzJd$ incorporates both the equivalence of the operators $(\upzd)_-$ and $(\upzJd)_-$ and the equivalence of the operators $(\upzd)_+$ and $(\upzJd)_+$. This leads to the following characterization of $q$--monogenicity depending on the location of the considered symplectic cell in the triangular decomposition of spinor space.

\begin{remark}
\label{equivalence-q-monog-pm}
We have that
\begin{itemize}
\item[(i)] in the inner part of the triangle, q--monogenicity is expressed by four operators: either $(\upz)_-$ or $(\upzJ)_-$, either $(\upz)_+$ or $(\upzJ)_+$, either $(\upzd)_-$ or $(\upzJd)_-$, and either  $(\upzd)_+$ or $(\upzJd)_+$; this makes sixteen possibilities;\\[-6mm]
\item[(ii)] on the left edge of the triangle ($\mathrm{Ker} \, P$), there are four possibilities to express q--monogenicity and to that end also four operators are needed: either $(\upz)_-$ or $(\upzJ)_-$, and $(\upz)_+$, either $(\upzd)_-$ or $(\upzJd)_-$, and $(\upzJd)_+$;\\[-6mm]
\item[(iii)] on the right edge of the triangle ($\mathrm{Ker} \, Q$), there are four possibilities to express q--monogenicity also needing four operators: either $(\upz)_-$ or $(\upzJ)_-$, and $(\upzJ)_+$, either $(\upzd)_-$ or $(\upzJd)_-$, and $(\upzd)_+$;\\[-6mm]
\item[(iv)] on the upper edge of the triangle only two operators are needed to express q--monogenicity and there are four possible choices: either $(\upz)_+$ or $(\upzJ)_+$, and either  $(\upzd)_+$ or $(\upzJd)_+$;\\[-6mm]
\item[(v)] at the lower vertex of the triangle also two operators are needed offering four possibilities: either $(\upz)_-$ or $(\upzJ)_-$, and either  $(\upzd)_-$ or $(\upzJd)_-$;\\[-6mm]
\item[(vi)] at the left upper vertex of the triangle two operators are needed leaving no choice at all: $(\upz)_+$ and $(\upzJd)_+$;\\[-6mm]
\item[(vii)] finally, at the right upper vertex, the two operators $(\upzJ)_+$ and $(\upzd)_+$ must be used.
\end{itemize}
\end{remark}


\section{Generalized gradients}
\label{gradients}


In this section we formulate  the equations of quaternionic Clifford analysis in an abstract way, following ideas of Stein and Weiss (\cite{stein}), and we establish a connection between these abstract equations and the actual equations of quaternionic Clifford analysis  considered in the preceding sections.
For the construction of these so--called {\em generalized gradients} we have to choose a fixed Euclidean vector space together with a hypercomplex structure, more explicitly, we fix:

\begin{itemize}
\item[(i)] the real vector space $V,$ of dimension $\dim V=4p,$ equipped with a Euclidean scalar product;\\[-7mm]
\item[(ii)] two anti--commuting complex structures $\mI$ and $\mJ$ on $V$ preserving this scalar product.
\end{itemize}

Introducing, in a natural way, a third complex structure by putting $\mK:=\mI \, \mJ$, we endow the vector space $V$ with a hypercomplex structure $\mQ=\{\mI,\mJ,\mK\}$. In this way we can reduce the symmetry group SO$(V)$ acting on $V$, to its subgroup
SO$_\mQ(V):=\{A \in SO(V) : A \mI = \mI A \ {\rm and} \ A \mJ = \mJ A\}$.
Let Sp$( p)$ denote the real Lie group of $p \times p$ matrices with  quaternion entries preserving the standard quaternionic inner product.
Then we have the following proposition.
 
\begin{proposition}\label{P1}
The group SO$_\mQ(V)$ is isomorphic to the Lie group Sp$( p)$.  
\end{proposition}

\pf
The following lemma shows that the complex structures $\mI_{4p}$ and $\mJ_{4p}$ used in the preceding sections can always be recovered by a suitable choice of an orthonormal basis in the vector space $V$.  Hence Proposition \ref{P1} follows from \cite{paper1}, Proposition 5.
\qed

\begin{lemma}
\label{basis}
There exists an orthonormal basis $\{e_i\}_{i=1}^{4p} $
of the vector space $V$ such that the complex structures $\mI$ and $\mJ$ are represented by the respective matrices $\mI_{4p}$ and $\mJ_{4p}$.
\end{lemma}

\pf
Let $B$ denote the Euclidean scalar product in the vector space $V$, and $B_c$ its complex bilinear extension to the complexification $V_c:=V\otimes\mC$. The maps $\mI$ and $\mJ$ also extend to complex linear maps on $V_c$; we keep the same notations for these extensions. The space $V_c$ decomposes as a direct sum $W \oplus W^\dagger$ of eigenspaces for $\mI$ with respective eigenvalues $\mp i.$ Both $W$ and $W^\dagger$ are isotropic with respect to $B_c$ and moreover $B_c$ induces a nondegenerate pairing between $W$ and $W^\dagger$. So if $\{\gf_j\}_{j=1}^{2p}$ is any basis of $W$, then there exists a basis $\{\gfd_k\}_{k=1}^{2p}$ for $W^\dagger$ such that 
$$
B_c(\gf_j,\gfd_k)=-\frac{1}{2}\delta_{jk}, \ j \neq k = 1,\ldots,2p
$$
Hence all elements in the union of both bases mutually anticommute; however mind the non--trivial relations $\gf_j \gfd_j+\gfd_j \gf_j = 1, j=1,\ldots,2p$.
Now the second complex structure $\mJ$ anti--commutes with $\mI$ and hence induces an isomorphism between $W$ and $W^\dagger$. Let us consider a complex bilinear anti--symmetric form $\omega$ on $W$ given by 
$$
\omega(\gf,\gf'):=2B_c(\gf,\mJ (\gf'))
$$
We can choose the basis
$\{\gf_j\}_{j=1}^{2p}$ for $W$ in such a way that $\omega(\gf_j,\gf_k)$ is trivial, except for the following cases:
$$
 \omega(\gf_{2j-1},\gf_{2j})=-\omega(\gf_{2j},
 \gf_{2j-1})=1,\,j=1,\ldots,p.
$$
Introducing a basis $\{e_j\}_{j=1}^{2p}$ by
$$
e_{2k-1} = \gfd_k - \gf_k \quad {\rm and} \quad e_{2k} = \frac{1}{i}(\gfd_k + \gf_k) \quad (k=1,\ldots,2p)
$$
it is easy to check that $\{e_j\}_{j=1}^{4p}$ is an orthonormal basis of $V$ and that the mappings $\mI$ and $\mJ$ have the required form.
\qed

\begin{lemma}
\label{V}
The vector space $V_c$ is a (complex) representation of {\rm SO}$_\mQ(V) \simeq$	 {\rm Sp}$( p)$.
As an {\rm Sp}$( p)$--module, $V_c$ decomposes as a sum of two copies of the defining representation $(1)_s$ of {\rm Sp}$( p)$.
\end{lemma}

\pf
As in the proof of Lemma \ref{basis}, the space $V_c$ decomposes as a direct sum $W \oplus W^\dagger$ of eigenspaces for $\mI$ with respective eigenvalues $\mp i.$
But, under the action of {\rm Sp}$(p)$, we have that $W^\dagger\simeq\mS^1_1$ and $W\simeq\mS^{2p-1}_1$ are both isomorphic to the defining representation $(1)_s$ of {\rm Sp}$( p)$.
\qed

\begin{remark}
For clarity's sake we note here  that the real symplectic Lie algebra $\gsp( p)$ of skew--symplectic $p\times p$--matrices with quaternion entries, is isomorphic with the so--called compact form $\gsp(2p,\mC) \cap \gu(2p)$ of the complex symplectic Lie algebra $\gsp(2p,\mC)$ (see e.g. \cite{paper1}, Proposition 6).
\end{remark}

\noindent
The construction of the Stein--Weiss gradients necessitates the projection on the components in the decomposition of a tensor product of two irreducible $\gsp(p)$-modules, which, in general, is a difficult problem. If one of the factors in the tensor product happens to be a {\it small} representation, it is possible to perform the decomposition explicitly using the Klimyk formula. In $\cite{SlSo}$, Section 5.9, it is shown how to achieve this if one of the factors is the defining representation of $\gsp(p).$ If one factor is a general representation with highest weight $\lambda,$ then the tensor product is multiplicity free and all summands have the form
$\lambda\pm\epsilon_i, i=1,\ldots,p,$ for which the result is dominant. Here $\epsilon_i$ has $1$ at the $i$--th place and zeros otherwise. In particular, we have the following lemma.
 
\begin{lemma}
\label{Klim}
Let $U_s$ denote an irreducible $\gsp(p)$-module with highest weight $(1,\ldots,1,0,\ldots,0)_s$ (with $s$ non--trivial entries). 
Let $V_c=W_1\oplus W_2$ be an irreducible decomposition with respect to $\gsp(p)$.
Then the tensor product  of $\gsp(p)$-modules $W_i\otimes U_s $ decomposes into irreducible components as
$$
U_{i,s+1}\oplus U_{i,s-1}\oplus U_{i,s}'
$$
where $U_{i,s\pm 1}\simeq U_{s\pm 1}$ and $U'_{i,s}$ is the Cartan product of $W_i$ and $U_s$. 

In what follows the invariant projections of $V_c\otimes U_s$ onto the summands $U_{i,s\pm 1}$  are denoted by $\pi_{i,\pm}$, respectively.
\end{lemma}

By the Stein--Weiss construction it is now possible to define four (invariant) first order differential operators
which are realized and used in quaternionic Clifford analysis.

\begin{definition}
\label{SWop}
Let the projections $\pi_{i,\pm}$ be defined as in Lemma \ref{Klim}.
Let $\Omega$ be an open domain in $V$. The first order differential operators $D_{i,\pm}: 
\mcC^\infty(\Omega, U_s) \mapsto \mcC^\infty(\Omega, U_{i,s\pm1})$
are defined by
\begin{equation}
\label{SWdef}
D_{i,\pm}(f) := \pi_{i,\pm}(\nabla(f))
\end{equation}
\end{definition}

It is clear that the abstract Stein--Weiss gradients $D_{i,\pm}$ defined above depend on the choice of the irreducible decomposition $V_c=W_1\oplus W_2$. 
But, fortunately, their common solutions do not depend on this choice. 

\begin{lemma}
\label{SWind} 
Let $D_{i,\pm}$, respectively \ $\tilde D_{i,\pm}$, be the differential operators defined as in Definition \ref{SWop} for the irreducible decomposition $V_c=W_1\oplus W_2$,  respectively \ $V_c=\tilde W_1\oplus\tilde W_2$. 
For each $f\in\mcC^\infty(\Omega, U_j)$, one has that $D_{1,\pm} (f)=D_{2,\pm} (f)=0$ if and only if $\tilde D_{1,\pm} (f)=\tilde D_{2,\pm} (f)=0$.
\end{lemma}

\pf
Let $\pi_{i,\pm}$, resp.\ $\tilde\pi_{i,\pm}$, be the projections defined as in Lemma \ref{Klim} for the irreducible decomposition $V_c=W_1\oplus W_2$,  resp.\ $V_c=\tilde W_1\oplus\tilde W_2$. 
Then $\pi_{\pm}=\pi_{1,\pm}+\pi_{2,\pm}$ are the projections of $V_c\otimes U_j$ onto the isotypic components $2U_{j\pm 1}$. Since the projections $\pi_{\pm}$
do not depend on the choice of irreducible decomposition of $V_c$ we have also $\pi_{\pm}=\tilde\pi_{1,\pm}+\tilde\pi_{2,\pm}$. Hence $\Ker (\pi_{1,\pm},\pi_{2,\pm})=\Ker (\tilde\pi_{1,\pm},\tilde\pi_{2,\pm})$. 
\qed

Now we show that quaternionic monogenicity can be expressed equivalently in terms of the Stein--Weiss gradients $D_{i,\pm}$. In this connection, let us recall that, by Proposition \ref{comp}, a spinor valued function is q--monogenic if and only if its symplectic components are.

\begin{theorem}
\label{main}
Let $D_{i,\pm}$ be the Stein--Weiss gradients defined as in Definition \ref{SWop} for $U_s=\mS^r_s$. 
Let $\Omega\subset\mR^{4p}$ be open. Then a~differentiable function $f:\Omega\to\mS^r_s$ is q--monogenic if and only if on $\Omega$ it satisfies the system
\begin{equation}
\label{SWeq}
D_{1,\pm} (f)=0,\ D_{2,\pm} (f)=0. 
\end{equation}  
\end{theorem} 

\pf ($a$) First assume that, for $\mS^r_s$-valued functions,  q--monogenicity can be characterized by the four operators $(\upz)_{\pm}$ and $(\upzd)_{\pm}$, see Remarks \ref{partial result} and \ref{equivalence-q-monog-pm} above. Without loss of generality, consider the case (i) of Remark \ref{equivalence-q-monog-pm}.
Furthermore, let $V_c=W \oplus W^\dagger$ be the decomposition into the eigenspaces for the complex structure $\mI$ with respective eigenvalues $\mp i$ and let $\tau_{i,\pm}$ be $\gsp(p)$-invariant isomorphisms of $U_{i,s\pm 1}$ onto $\mS^{r+2i-3}_{s\pm 1}$. 
Then, for $W_1=W$ and $W_2=W^\dagger$, we have that $(\upzd)_{\pm}$ and $(\upz)_{\pm}$  are respectively equal to $\tau_{1,\pm}\circ D_{1,\pm}$ and $\tau_{2,\pm}\circ D_{2,\pm}$ up to non-zero multiples, see \cite[Section 4]{partII}.  
Hence, in this case, the required equivalence is proven.

\smallskip\noindent
($b$) Now consider, say, the case (ii) of Remark \ref{equivalence-q-monog-pm}, i.e., $r=s<p$. 
Let $f:\Omega\to\mS^r_r$ be a~differentiable function. 
By Lemma \ref{PQF}, the function $f$ is q--monogenic if and only if $Qf$ is. By construction it is clear that the system $D_{1,\pm} (f)=0=D_{2,\pm} (f)$ is equivalent to the system $\tilde D_{1,\pm} (Qf)=0=\tilde D_{2,\pm} (Qf)$, where 
$\tilde D_{i,\pm}$ are the Stein--Weiss gradients defined for $\tilde U_s=\mS^{r+2}_r$.  
We complete the proof by applying ($a$) to $Qf$. The remaining cases can be proven similarly.    
\qed


\section{Characterization of q--monogenic functions for p = 2}
\label{lowdimension}


In the last section we illustrate our findings above by the explicit calculations of the systems of equations corresponding with q--monogenic functions defined in $\mR^8$ and taking values in spinor space $\mS$ with $p=2$. The  triangular scheme for the decomposition of this spinor space looks as follows:

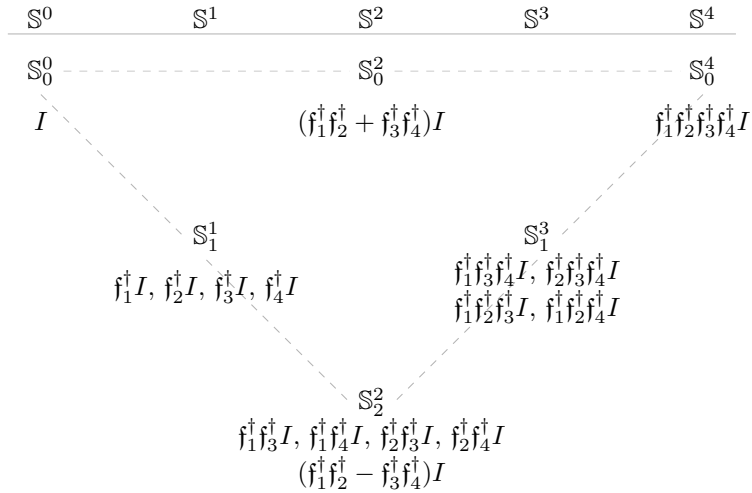
\begin{center}
\begin{tikzpicture}[scale=1.1]
\node at (0,0) {$\mS^0$};
\node at (2,0) {$\mS^1$}; 
\node at (4,0) {$\mS^2$};
\node at (6,0) {$\mS^3$}; 
\node at (8,0) {$\mS^4$};
\draw[draw opacity = 0.3] (-0.4,-0.2) -- (8.4,-0.2);
\node at (0,-0.65) (S00) {$\mS_0^0$};
\node at (0,-1.25) (I) {$I$};
\node at (4,-0.65) (S02) {$\mS_0^2$};
\node at (4,-1.25) (QI) {$(\gfd_1 \gfd_2 + \gfd_3 \gfd_4) I$};
\node at (8,-0.65) (S04) {$\mS_0^4$};
\node at (8,-1.25) (QQI) {$\gfd_1 \gfd_2 \gfd_3 \gfd_4 I$};
\draw[dashed,draw opacity = 0.2] (S00.east) to (S02.west);
\draw[dashed,draw opacity = 0.2] (S02.east) to (S04.west);
\node at (2,-2.65) (S11) {$\mS^1_1$}; 
\node at (6,-2.65) (S13) {$\mS^3_1$}; 
\draw[dashed,draw opacity = 0.3] (S00.south) to (S11.west);
\draw[dashed,draw opacity = 0.3] (S04.south) to (S13.east);
\node at (2,-3.25) {$\gfd_1 I$, $\gfd_2 I$, $\gfd_3 I$, $\gfd_4 I$};
\node at (6,-3.05) {$\gfd_1 \gfd_3 \gfd_4 I$, $\gfd_2 \gfd_3 \gfd_4 I$};
\node at (6,-3.5) {$\gfd_1 \gfd_2 \gfd_3 I$, $\gfd_1 \gfd_2 \gfd_4 I$};
\node at (4,-4.65) (S22) {$\mS^2_2$}; 
\draw[dashed,draw opacity = 0.3] (S11.south) to (S22.west);
\draw[dashed,draw opacity = 0.3] (S13.south) to (S22.east);
\node at (4,-5.05)  {$\gfd_1 \gfd_3 I$, $\gfd_1 \gfd_4 I$, $\gfd_2 \gfd_3 I$, $\gfd_2 \gfd_4 I$};
\node at (4,-5.5) {$(\gfd_1 \gfd_2 - \gfd_3 \gfd_4) I$};
\end{tikzpicture}
\end{center}

\noindent
{\bf Case A:} The function $F_0^0 : \mR^8 \longrightarrow \mS_0^0$ has the form $F_0^0 = \phi \,  I$, whence we obtain
\begin{itemize}

\item $\upz F_0^0 = \left( \gfd_1 \, \p_{z_1}\phi + \gfd_2 \, \p_{z_2}\phi + \gfd_3 \, \p_{z_3}\phi + \gfd_4 \, \p_{z_4}\phi \right) I$\\[2mm]
which splits into\\[2mm]
$(\upz)_+ F_0^0 =   \left( \gfd_1 \, \p_{z_1}\phi + \gfd_2 \, \p_{z_2}\phi + \gfd_3 \, \p_{z_3}\phi + \gfd_4 \, \p_{z_4}\phi \right) I$\\[1mm]
$(\upzd)_- F_0^0 = 0$
\end{itemize}
and
\begin{itemize}
\item $\upzJd F_0^0 = \left(  \gf_1 \, \p_{\olz_2}\phi -  \gfd_2 \, \p_{\olz_1}\phi + \gfd_3 \, \p_{\olz_4}\phi - \gfd_4 \, \p_{\olz_3}\phi \right) I$\\[2mm]
which splits into\\[2mm]
$(\upzJd)_+ F_0^0 = \left( \gf_1 \, \p_{\olz_2}\phi -  \gfd_2 \, \p_{\olz_1}\phi + \gfd_3 \, \p_{\olz_4}\phi - \gfd_4 \, \p_{\olz_3}\phi \right)  I$\\[1mm]
$(\upzJd)_- F_0^0 = 0$
\end{itemize}
while, trivially
\begin{itemize}
\item $\upzd F_0^0 = 0 = \upzJ F_0^0 = 0$ 
\end{itemize}
It follows that
\begin{itemize}
\item $\upz F_0^0 = 0 \Longleftrightarrow \p_{z_1}\phi = 0 ,   \p_{z_2}\phi = 0 , \p_{z_3}\phi = 0 , \p_{z_4}\phi = 0$
\item $\upzJd F_0^0 =  0 \Longleftrightarrow  \p_{\olz_2}\phi = 0 ,  \p_{\olz_1}\phi = 0 , \p_{\olz_3}\phi = 0 , \p_{\olz_4}\phi = 0$
\end{itemize}
In other words: $F_0^0$ is q--monogenic if and only if $\phi$ is a constant function. Moreover, this observation is valid regardless the dimension. Seen the fact that all symplectic cells on the same row in the spinor triangle are realizations of the same Sp$( p)$ representation, we now expect the q--monogenic functions with values in $\mS^2_0$ (case C) and $\mS^4_0$ (case F) to be constant functions too. We also see that for expressing the q--monogenicity of $F_0^0$ the two operators $(\upz)_+$ and $(\upzJd)_+$ need to be used (left upper vertex of the triangle).\\[-2mm]

\noindent
{\bf Case B:} The function $F_1^1 : \mR^8 \longrightarrow \mS_1^1$ has the form 
$$
F_1^1 = (\phi_1 \,  \gfd_1 +  \phi_2 \,  \gfd_2 + \phi_3 \,  \gfd_3 + \phi_4 \,  \gfd_4 ) \, I
$$
so there holds
\begin{itemize}

\item $\upz F_1^1$\\[2mm]
$ = \left(  \gfd_1 \gfd_3 \, (\p_{z_1}\phi_3 - \p_{z_3}\phi_1) + \gfd_1 \gfd_4 \, (\p_{z_1}\phi_4 - \p_{z_4}\phi_1)+ \gfd_2 \gfd_3 \, (\p_{z_2}\phi_3 - \p_{z_3}\phi_2) + \gfd_2 \gfd_4 \, (\p_{z_2}\phi_4 - \p_{z_4}\phi_2) \right) I \\
\hspace*{3mm}+ \left(  \gfd_1 \gfd_2 ( \p_{z_1}\phi_2 - \p_{z_2}\phi_1) + \gfd_3 \gfd_4 ( \p_{z_3}\phi_4 - \p_{z_4}\phi_3) \right) I$\\[2mm]
which splits into\\[2mm]
$(\upz)_+ F_1^1$\\[2mm]
$ =  \left(  \gfd_1 \gfd_3 \, (\p_{z_1}\phi_3 - \p_{z_3}\phi_1) + \gfd_1 \gfd_4 \, (\p_{z_1}\phi_4 - \p_{z_4}\phi_1)+ \gfd_2 \gfd_3 \, (\p_{z_2}\phi_3 - \p_{z_3}\phi_2) + \gfd_2 \gfd_4 \, (\p_{z_2}\phi_4 - \p_{z_4}\phi_2) \right) I$ \\
\hspace*{3mm}$ + \left( \onehalf  (\gfd_1 \gfd_2 - \gfd_3 \gfd_4)( \p_{z_1}\phi_2 - \p_{z_2}\phi_1) - \onehalf(\gfd_1 \gfd_2 - \gfd_3 \gfd_4) ( \p_{z_3}\phi_4 - \p_{z_4}\phi_3) \right) I$\\[1mm]
$(\upz)_- F_1^1 = \left( \onehalf  (\gfd_1 \gfd_2 + \gfd_3 \gfd_4)( \p_{z_1}\phi_2 - \p_{z_2}\phi_1) + \onehalf(\gfd_1 \gfd_2 + \gfd_3 \gfd_4) ( \p_{z_3}\phi_4 - \p_{z_4}\phi_3) \right) I$

\item $\upzd F_1^1 =  \left( \p_{\olz_1} \phi_1 + \p_{\olz_2} \phi_2 + \p_{\olz_3} \phi_3 + \p_{\olz_4} \phi_4 \right) I$\\[2mm]
which splits into\\[2mm]
$(\upzd)_+ F_1^1 = 0$\\[1mm]
$(\upzd)_- F_1^1 = \left( \p_{\olz_1} \phi_1 + \p_{\olz_2} \phi_2 + \p_{\olz_3} \phi_3 + \p_{\olz_4} \phi_4 \right) I$

\item $\upzJ F_1^1 = \left( \p_{z_2} \phi_1 - \p_{z_1} \phi_2 + \p_{z_4} \phi_3 - \p_{z_3} \phi_4 \right) I$\\[2mm]
which splits into\\[2mm]
$(\upzJ)_+ F_1^1 = 0$\\[1mm]
$(\upzJ)_- F_1^1 = \left( \p_{z_2} \phi_1 - \p_{z_1} \phi_2 + \p_{z_4} \phi_3 - \p_{z_3} \phi_4 \right) I$
\end{itemize}
and
\begin{itemize}
\item $\upzJd F_1^1$\\[2mm]
$ =  \left(  \gfd_1 \gfd_3 \, (\p_{\olz_2}\phi_3 - \p_{\olz_4}\phi_1) + \gfd_1 \gfd_4 \, (\p_{\olz_2}\phi_4 +\p_{\olz_3} \phi_1 ) + \gfd_2 \gfd_3 \, (- \p_{\olz_1}\phi_3 - \p_{\olz_4}\phi_2) + \gfd_2 \gfd_4 \, (\p_{\olz_3}\phi_2 - \p_{\olz_1}\phi_4) \right) I $\\
\hspace*{3mm} $ + \left(  \gfd_1 \gfd_2 ( \p_{\olz_2}\phi_2 + \p_{\olz_1}\phi_1) + \gfd_3 \gfd_4 ( \p_{\olz_4}\phi_4 + \p_{\olz_3}\phi_3) \right) I$\\[2mm]
which splits into\\[2mm]
$(\upzJd)_+ F_1^1$\\[2mm]
$ =  \left(  \gfd_1 \gfd_3 \, (\p_{\olz_2}\phi_3 - \p_{\olz_4}\phi_1) + \gfd_1 \gfd_4 \, (\p_{\olz_2}\phi_4 +\p_{\olz_3} \phi_1 ) + \gfd_2 \gfd_3 \, (- \p_{\olz_1}\phi_3 - \p_{\olz_4}\phi_2) + \gfd_2 \gfd_4 \, (\p_{\olz_3}\phi_2 - \p_{\olz_1}\phi_4) \right) I $\\
\hspace*{3mm} $ + \left(  \onehalf(\gfd_1 \gfd_2 - \gfd_3 \gfd_4) ( \p_{\olz_2}\phi_2 + \p_{\olz_1}\phi_1) - \onehalf (\gfd_1 \gfd_2 - \gfd_3 \gfd_4) ( \p_{\olz_4}\phi_4 + \p_{\olz_3}\phi_3) \right) I$\\[1mm]
$(\upzJd)_- F_1^1 = \left(  \onehalf(\gfd_1 \gfd_2 + \gfd_3 \gfd_4) ( \p_{\olz_2}\phi_2 + \p_{\olz_1}\phi_1) + \onehalf (\gfd_1 \gfd_2 + \gfd_3 \gfd_4) ( \p_{\olz_4}\phi_4 + \p_{\olz_3}\phi_3) \right) I$

\end{itemize}
It follows that
\begin{itemize}

\item $\upz F_1^1 = 0 \Longleftrightarrow \p_{z_1}\phi_3 - \p_{z_3}\phi_1 = 0 \ , \ \p_{z_1}\phi_4 - \p_{z_4}\phi_1= 0 \ , \ \p_{z_2}\phi_3 - \p_{z_3}\phi_2 = 0 \ , \ \p_{z_2}\phi_4 - \p_{z_4}\phi_2 = 0 \ , \  \p_{z_1}\phi_2 - \p_{z_2}\phi_1 = 0 \ , \
 \p_{z_3}\phi_4 - \p_{z_4}\phi_3 = 0$
 
\item $ \upzd F_1^1 = 0 \Longleftrightarrow \p_{\olz_1} \phi_1 + \p_{\olz_2} \phi_2 + \p_{\olz_3} \phi_3 + \p_{\olz_4} \phi_4 = 0$

\item $ \upzJ F_1^1 = 0 \Longleftrightarrow \p_{z_2} \phi_1 - \p_{z_1} \phi_2 + \p_{z_4} \phi_3 - \p_{z_3} \phi_4 = 0$
 
\item $\upzJd F_1^1 =  0 \Longleftrightarrow  \p_{\olz_2}\phi_3 - \p_{\olz_4}\phi_1 = 0 \ , \  \p_{\olz_2}\phi_4 +\p_{\olz_3} \phi_1 = 0 \ , \  \p_{\olz_1}\phi_3 + \p_{\olz_4}\phi_2 = 0 \ , \  \p_{\olz_3}\phi_2 - \p_{\olz_1}\phi_4 = 0 \ , \
 \p_{\olz_2}\phi_2 + \p_{\olz_1}\phi_1 = 0 \ , \   \p_{\olz_4}\phi_4 + \p_{\olz_3}\phi_3 = 0 $
 
\end{itemize}
We see that $\upz F_1^1 = 0 $ implies $\upzJ F_1^1 = 0 $, but not conversely, and that  $\upzJd F_1^1 = 0 $ implies $\upzd F_1^1 = 0 $, but also not conversely.
We also see that four operators are needed to express q--monogenicty: $(\upz)_+$, $(\upzJd)_+$, either $(\upz)_-$ or $(\upzJ)_-$ which both lead to the same system of equations, and either $(\upzd)_-$ or $(\upzJd)_-$ which both lead also to the same system of equations as well. Moreover we expect q--monogenic functions with values in $\mS^3_1$ (case E) to satisfy the same system of equations as $F^1_1$ does.\\[-2mm]

\noindent
{\bf Case C:} The function $F_0^2 : \mR^8 \longrightarrow \mS_0^2$ has the form 
$$
F_0^2 = \phi \, \left( \gfd_1 \gfd_2 + \gfd_3 \gfd_4 \right) I
$$
whence there holds
\begin{itemize}
\item $\upz F_0^2 = \left( \gfd_1 \gfd_3 \gfd_4 \, \p_{z_1}\phi + \gfd_2 \gfd_3 \gfd_4 \, \p_{z_2}\phi + \gfd_1 \gfd_2 \gfd_3 \, \p_{z_3}\phi + \gfd_1 \gfd_2 \gfd_4 \, \p_{z_4}\phi \right) I$\\[2mm]
which splits into\\[2mm]
$(\upz)_+ F_0^2 =   \left( \gfd_1 \gfd_3 \gfd_4 \, \p_{z_1}\phi + \gfd_2 \gfd_3 \gfd_4 \, \p_{z_2}\phi + \gfd_1 \gfd_2 \gfd_3 \, \p_{z_3}\phi + \gfd_1 \gfd_2 \gfd_4 \, \p_{z_4}\phi \right) I$\\[1mm]
$(\upzd)_- F_0^2 = 0$

\item $\upzd F_0^2 =  \left( \gfd_2 \, \p_{\olz_1}\phi + \gfd_1 \, \p_{\olz_2}\phi + \gfd_4 \, \p_{\olz_3}\phi + \gfd_3 \, \p_{\olz_4}\phi \right) I$\\[2mm]
which splits into\\[2mm]
$(\upzd)_+ F_0^2 = \left( \gfd_2 \, \p_{\olz_1}\phi + \gfd_1 \, \p_{\olz_2}\phi + \gfd_4 \, \p_{\olz_3}\phi + \gfd_3 \, \p_{\olz_4}\phi \right) I$\\[1mm]
$(\upzd)_- F_0^2 = 0$

\item $\upzJ F_0^2 = \left(  \gfd_1 \, \p_{z_1}\phi + \gfd_2 \, \p_{z_2}\phi  + \gfd_3 \, \p_{z_3}\phi + \gfd_4 \, \p_{z_4}\phi  \right) I$\\[2mm]
which splits into\\[2mm]
$(\upzJ)_+ F_0^2 = \left(  \gfd_1 \, \p_{z_1}\phi + \gfd_2 \, \p_{z_2}\phi  + \gfd_3 \, \p_{z_3}\phi + \gfd_4 \, \p_{z_4}\phi  \right) I$\\[1mm]
$(\upzJ)_- F_0^2 = 0 $
\end{itemize}
and
\begin{itemize}
\item $\upzJd F_0^2 = \left( \gfd_1 \gfd_3 \gfd_4 \, \p_{\olz_2}\phi - \gfd_2 \gfd_3 \gfd_4 \, \p_{\olz_1}\phi + \gfd_1 \gfd_2 \gfd_3 \, \p_{\olz_4}\phi - \gfd_1 \gfd_2 \gfd_4 \, \p_{\olz_3}\phi \right) I$\\[2mm]
which splits into\\[2mm]
$(\upzJd)_+ F_0^2 =  \left( \gfd_1 \gfd_3 \gfd_4 \, \p_{\olz_2}\phi - \gfd_2 \gfd_3 \gfd_4 \, \p_{\olz_1}\phi + \gfd_1 \gfd_2 \gfd_3 \, \p_{\olz_4}\phi - \gfd_1 \gfd_2 \gfd_4 \, \p_{\olz_3}\phi \right) I$\\[1mm]
$(\upzJd)_- F_0^2 = 0$
\end{itemize}
It follows that
\begin{itemize}
\item $\upz F_0^2 = 0 \Longleftrightarrow \p_{z_1}\phi = 0 ,   \p_{z_2}\phi = 0 , \p_{z_3}\phi = 0 , \p_{z_4}\phi = 0$
\item $\upzd F_0^2 = 0 \Longleftrightarrow \p_{\olz_1}\phi = 0 ,  \p_{\olz_2}\phi = 0 , \p_{\olz_3}\phi = 0 , \p_{\olz_4}\phi$
\item $\upzJ F_0^2 = 0 \Longleftrightarrow \p_{z_1}\phi = 0 ,   \p_{z_2}\phi = 0 , \p_{z_3}\phi = 0 , \p_{z_4}\phi = 0$
\item $\upzJd F_0^2 =  0 \Longleftrightarrow  \p_{\olz_2}\phi = 0 ,  \p_{\olz_1}\phi = 0 , \p_{\olz_4}\phi = 0 , \p_{\olz_3}\phi = 0$
\end{itemize}
in other words: $F_0^2$ is q--monogenic if and only if $\phi$ is a constant function. Besides, this property remains valid for all functions with values in a (one--dimensional) cell on the upper edge of the spinor space triangle, and this regardless the dimension: if $F^r_0 : \mR^{4p} \longrightarrow \mS^r_0$ is q--monogenic, then it reduces to a constant function. We see that the operators $\upz$ and $\upzJ$ on the one hand, and the operators $\upzd$ and $\upzJd$ on the other hand, lead to the same system of equations. We also see that only two operators are needed to express q--monogenicity: either $(\upz)_+$ or $(\upzJ)_+$, both leading to the same system of equations, and either  $(\upzd)_+$ or $(\upzJd)_+$, both leading also to the same system of equations.\\[-2mm]

\noindent
{\bf Case D:} The function $F_2^2 : \mR^8 \longrightarrow \mS_2^2$ has the form 
$$
F_2^2 = \left(\phi_{13} \,  \gfd_1 \gfd_3 +  \phi_{14 }\,  \gfd_1  \gfd_4 + \phi_{23} \,  \gfd_2 \gfd_3 + \phi_{24} \, \gfd_2 \gfd_4  + \phi  \left( \gfd_1 \gfd_2 - \gfd_3 \gfd_4\right)\right)  I
$$
whence there holds
\begin{itemize}
\item $\upz F_2^2 =  (  \gfd_1 \gfd_2 \gfd_3 \, (- \p_{z_2}\phi_{13} + \p_{z_1}\phi_{23} + \p_{z_3} \phi) + \gfd_1 \gfd_2 \gfd_4 \, (\p_{z_1}\phi_{24} - \p_{z_2} \phi_{14} + \p_{z_4} \phi) \\
\hspace*{14mm} + \gfd_1 \gfd_3 \gfd_4 \, (- \p_{z_3}\phi_{14} - \p_{z_1}\phi + \p_{z_4} \phi_{13}) + \gfd_2 \gfd_3 \gfd_4 \, (\p_{z_4} \phi_{23} - \p_{z_3}\phi_{24} - \p_{z_2}\phi)) I$\\[2mm]
which splits into\\[2mm]
$(\upz)_+ F_2^2 =  0$\\[1mm]
$(\upzd)_- F_2^2 = (  \gfd_1 \gfd_2 \gfd_3 \, (- \p_{z_2}\phi_{13} + \p_{z_1}\phi_{23} + \p_{z_3} \phi) + \gfd_1 \gfd_2 \gfd_4 \, (\p_{z_1}\phi_{24} - \p_{z_2} \phi_{14} + \p_{z_4} \phi) \\
\hspace*{19mm} + \gfd_1 \gfd_3 \gfd_4 \, (- \p_{z_3}\phi_{14} - \p_{z_1}\phi + \p_{z_4} \phi_{13}) + \gfd_2 \gfd_3 \gfd_4 \, (\p_{z_4} \phi_{23} - \p_{z_3}\phi_{24} - \p_{z_2}\phi)) I$

\item $\upzd F_2^2 =  ( \gfd_3 \, ( \p_{\olz_1}\phi_{13} + \p_{\olz_2}\phi_{23} + \p_{\olz_4} \phi) + \gfd_2 \, (- \p_{\olz_4}\phi_{24} - \p_{\olz_3} \phi_{23} + \p_{\olz_1} \phi) \\
\hspace*{14mm} + \gfd_4 \, ( \p_{\olz_1}\phi_{14} + \p_{\olz_2}\phi_{24} - \p_{\olz_3} \phi) + \gfd_1 \, (- \p_{\olz_3} \phi_{13} - \p_{\olz_4}\phi_{14} - \p_{\olz_2}\phi)) I$\\[2mm]
which splits into\\[2mm]
$(\upzd)_+ F_2^2 = 0$\\[1mm]
$(\upzd)_- F_2^2 =  ( \gfd_3 \, ( \p_{\olz_1}\phi_{13} + \p_{\olz_2}\phi_{23} + \p_{\olz_4} \phi) + \gfd_2 \, (- \p_{\olz_4}\phi_{24} - \p_{\olz_3} \phi_{23} + \p_{\olz_1} \phi) \\
\hspace*{19mm} + \gfd_4 \, ( \p_{\olz_1}\phi_{14} + \p_{\olz_2}\phi_{24} - \p_{\olz_3} \phi) + \gfd_1 \, (- \p_{\olz_3} \phi_{13} - \p_{\olz_4}\phi_{18} - \p_{\olz_2}\phi)) I$

\item $\upzJ F_2^2 = ( \gfd_3 \, ( \p_{z_2}\phi_{13} - \p_{z_1}\phi_{23} - \p_{z_3} \phi) + \gfd_2 \, ( \p_{z_3}\phi_{24} - \p_{z_4} \phi_{23} + \p_{z_2} \phi) \\
\hspace*{15mm} + \gfd_4 \, ( \p_{z_2}\phi_{14} - \p_{z_1}\phi_{24} - \p_{z_4} \phi) + \gfd_1 \, ( \p_{z_3} \phi_{14} - \p_{z_4}\phi_{13} + \p_{z_1}\phi)) I$\\[2mm]
which splits into\\[2mm]
$(\upzJ)_+ F_2^2 = 0$\\[1mm]
$(\upzJ)_- F_2^2 = ( \gfd_3 \, ( \p_{z_2}\phi_{13} - \p_{z_1}\phi_{23} - \p_{z_3} \phi) + \gfd_2 \, ( \p_{z_3}\phi_{24} - \p_{z_4} \phi_{23} + \p_{z_2} \phi) \\
\hspace*{19mm} + \gfd_4 \, ( \p_{z_2}\phi_{14} - \p_{z_1}\phi_{24} - \p_{z_4} \phi) + \gfd_1 \, ( \p_{z_3} \phi_{14} - \p_{z_4}\phi_{13} + \p_{z_1}\phi)) I$
\end{itemize}
and
\begin{itemize}
\item $\upzJd F_2^2 =  (  \gfd_1 \gfd_2 \gfd_3 \, ( \p_{\olz_1}\phi_{13} + \p_{\olz_2}\phi_{23} + \p_{\olz_4} \phi) + \gfd_1 \gfd_2 \gfd_4 \, (\p_{\olz_2}\phi_{24} + \p_{\olz_1} \phi_{14} - \p_{\olz_3} \phi) \\
\hspace*{16mm} + \gfd_1 \gfd_3 \gfd_4 \, (- \p_{\olz_4}\phi_{14} - \p_{\olz_2}\phi - \p_{\olz_3} \phi_{13}) + \gfd_2 \gfd_3 \gfd_4 \, (- \p_{\olz_3} \phi_{23} - \p_{\olz_4}\phi_{24} + \p_{\olz_1}\phi)) I$\\[2mm]
which splits into\\[2mm]
$(\upzJd)_+ F_2^2 =  0$\\[1mm]
$(\upzJd)_- F_2^2 = (  \gfd_1 \gfd_2 \gfd_3 \, ( \p_{\olz_1}\phi_{13} + \p_{\olz_2}\phi_{23} + \p_{\olz_4} \phi) + \gfd_1 \gfd_2 \gfd_4 \, (\p_{\olz_2}\phi_{24} + \p_{\olz_1} \phi_{14} - \p_{\olz_3} \phi) \\
\hspace*{20mm} + \gfd_1 \gfd_3 \gfd_4 \, (- \p_{\olz_4}\phi_{14} - \p_{\olz_2}\phi - \p_{\olz_3} \phi_{13}) + \gfd_2 \gfd_3 \gfd_4 \, (- \p_{\olz_3} \phi_{23} - \p_{\olz_4}\phi_{24} + \p_{\olz_1}\phi)) I $
\end{itemize}
It follows that
\begin{itemize}

\item $\upz F_2^2 = 0 \Longleftrightarrow - \p_{z_2}\phi_{13} + \p_{z_1}\phi_{23} + \p_{z_3} \phi = 0 \ , \  \p_{z_1}\phi_{24} - \p_{z_2} \phi_{14} + \p_{z_4} \phi = 0 \ , \\
\hspace*{22mm} - \p_{z_3}\phi_{14} - \p_{z_1}\phi + \p_{z_4} \phi_{13} = 0 \ , \ \p_{z_4} \phi_{23} - \p_{z_3}\phi_{24} - \p_{z_2}\phi = 0$
 
\item $ \upzd F_2^2 = 0 \Longleftrightarrow  \p_{\olz_1}\phi_{13} + \p_{\olz_2}\phi_{23} + \p_{\olz_4} \phi = 0 \ , \  - \p_{\olz_4}\phi_{24} - \p_{\olz_3} \phi_{23} + \p_{\olz_1} \phi = 0 \ , \\
\hspace*{23mm} \p_{\olz_1}\phi_{14} + \p_{\olz_2}\phi_{24} - \p_{\olz_3} \phi = 0 \ , \  - \p_{\olz_3} \phi_{13} - \p_{\olz_4}\phi_{14} - \p_{\olz_2}\phi = 0$ 
 
\item $ \upzJ F_2^2 = 0 \Longleftrightarrow  \p_{z_2}\phi_{13} - \p_{z_1}\phi_{23} - \p_{z_3} \phi = 0 \ , \  \p_{z_3}\phi_{24} - \p_{z_4} \phi_{23} + \p_{z_2} \phi = 0 \ , \\
\hspace*{23mm}  \p_{z_2}\phi_{14} - \p_{z_1}\phi_{24} - \p_{z_4} \phi = 0 \ , \  \p_{z_3} \phi_{14} - \p_{z_4}\phi_{13} + \p_{z_1}\phi = 0$
 
\item $\upzJd F_2^2 =  0 \Longleftrightarrow  \p_{\olz_1}\phi_{13} + \p_{\olz_2}\phi_{23} + \p_{\olz_4} \phi = 0 \ , \ \p_{\olz_2}\phi_{24} + \p_{\olz_1} \phi_{14} - \p_{\olz_3} \phi = 0 \ ,\\
\hspace*{23mm} - \p_{\olz_4}\phi_{14} - \p_{\olz_2}\phi - \p_{\olz_3} \phi_{13} = 0 \ , \  - \p_{\olz_3} \phi_{23} - \p_{\olz_4}\phi_{24} + \p_{\olz_1}\phi = 0$
 
\end{itemize}
We see that the operators $\upz$ and $\upzJ$ on the one hand, and the operators $\upzd$ and $\upzJd$ on the other hand lead to the same system of equations.
We also see that two operators are needed to express q--monogenicity: either $(\upz)_-$ or $(\upzJ)_-$, both leading to the same system of equations, and either  $(\upzd)_-$ or $(\upzJd)_-$, both leading also to the same system of equations.\\[-2mm]

\noindent
{\bf Case E:} The function $F_1^3 : \mR^8 \longrightarrow \mS_1^3$ has the form 
$$
F_1^3 = \left(\phi_1 \,  \gfd_1 \gfd_3 \gfd_4 +  \phi_2 \,  \gfd_2  \gfd_3 \gfd_4 + \phi_3 \,  \gfd_1 \gfd_2 \gfd_3 + \phi_4 \,  \gfd_1 \gfd_2 \gfd_4 \right)  I
$$
whence there holds
\begin{itemize}

\item $\upz F_1^3 = \left(  \p_{z_1}\phi_2 - \p_{z_2}\phi_1 + \p_{z_3}\phi_4 - \p_{z_4}\phi_3) \right) \gfd_1 \gfd_2 \gfd_3 \gfd_4 I $\\[2mm]
which splits into\\[2mm]
$(\upz)_+ F_1^3 =  0$\\[1mm]
$(\upz)_- F_1^3 = \left(  \p_{z_1}\phi_2 - \p_{z_2}\phi_1 + \p_{z_3}\phi_4 - \p_{z_4}\phi_3) \right) \gfd_1 \gfd_2 \gfd_3 \gfd_4 I $

\item $\upzd F_1^3 =  (  \gfd_1 \gfd_3 \, (\p_{\olz_4}\phi_1 - \p_{\olz_2}\phi_3) + \gfd_1 \gfd_4 \, (-\p_{\olz_2}\phi_4 - \p_{\olz_3} \phi_3 ) + \gfd_2 \gfd_3 \, ( \p_{\olz_1}\phi_3 + \p_{\olz_4}\phi_2))I$\\
\hspace*{12mm} $ + (\gfd_2 \gfd_4 \, (- \p_{\olz_3}\phi_2 + \p_{\olz_1}\phi_4)  +  \gfd_1 \gfd_2 ( \p_{\olz_3}\phi_3 + \p_{\olz_4}\phi_4) + \gfd_3 \gfd_4 ( \p_{\olz_1}\phi_1 + \p_{\olz_2}\phi_2) ) I$\\[2mm]
which splits into\\[2mm]
$(\upzd)_+ F_1^3 = (  \gfd_1 \gfd_3 \, (\p_{\olz_4}\phi_1 - \p_{\olz_2}\phi_3) + \gfd_1 \gfd_4 \, (-\p_{\olz_2}\phi_4 - \p_{\olz_3} \phi_1 ))I$\\
\hspace*{17mm} $+ (\gfd_2 \gfd_3 \, ( \p_{\olz_1}\phi_3 + \p_{\olz_4}\phi_2))I + \gfd_2 \gfd_4 \, (- \p_{\olz_3}\phi_2 + \p_{\olz_1}\phi_4) ) I$ \\
\hspace*{17mm} $+ \left( \onehalf ( \gfd_1 \gfd_2 - \gfd_3 \gfd_4) ( \p_{\olz_3}\phi_3 + \p_{\olz_4}\phi_4) - \onehalf ( \gfd_1 \gfd_2 - \gfd_3 \gfd_4) ( \p_{\olz_1}\phi_1 + \p_{\olz_2}\phi_2) \right) I$\\[1mm]
$(\upzd)_- F_1^3 = \left( \onehalf ( \gfd_1 \gfd_2 + \gfd_3 \gfd_4) ( \p_{\olz_3}\phi_3 + \p_{\olz_4}\phi_4) + \onehalf ( \gfd_1 \gfd_2 + \gfd_3 \gfd_4) ( \p_{\olz_1}\phi_1 + \p_{\olz_2}\phi_2) \right) I$

\item $\upzJ F_1^3 = (  \gfd_1 \gfd_3 \, (- \p_{z_3}\phi_1 + \p_{z_1}\phi_3) + \gfd_1 \gfd_4 \, (\p_{z_1}\phi_4 - \p_{z_4} \phi_1 ) + \gfd_2 \gfd_3 \, ( \p_{z_2}\phi_3 - \p_{z_3}\phi_2) ) I$\\
\hspace*{12mm} $ + (\gfd_2 \gfd_4 \, (- \p_{z_4}\phi_2 + \p_{z_2}\phi_4) +  \gfd_1 \gfd_2 ( \p_{z_4}\phi_3 - \p_{z_3}\phi_4) + \gfd_3 \gfd_4 ( \p_{z_2}\phi_1 - \p_{z_1}\phi_2) ) I$\\[2mm]
which splits into\\[2mm]
$(\upzJ)_+ F_1^3 = (  \gfd_1 \gfd_3 \, (- \p_{z_3}\phi_1 + \p_{z_1}\phi_3) + \gfd_1 \gfd_4 \, (\p_{z_1}\phi_4 - \p_{z_4} \phi_1 ))I$\\
\hspace*{18mm} $+ \gfd_2 \gfd_3 \, ( \p_{z_2}\phi_3 - \p_{z_3}\phi_2) + \gfd_2 \gfd_4 \, (- \p_{z_4}\phi_2 + \p_{z_2}\phi_4) ) I$ \\
\hspace*{18mm} $ + \left( \onehalf( \gfd_1 \gfd_2 - \gfd_3 \gfd_4)( \p_{z_4}\phi_3 - \p_{z_3}\phi_4) - \onehalf (\gfd_1 \gfd_2 - \gfd_3 \gfd_4) ( \p_{z_2}\phi_1 - \p_{z_1}\phi_2) \right) I$\\[1mm]
$(\upzJ)_- F_1^3 = \left( \onehalf( \gfd_1 \gfd_2 + \gfd_3 \gfd_4)( \p_{z_4}\phi_3 - \p_{z_3}\phi_4) + \onehalf (\gfd_1 \gfd_2 + \gfd_3 \gfd_4) ( \p_{z_2}\phi_1 - \p_{z_1}\phi_2) \right) I$
\end{itemize}
and
\begin{itemize}
\item $\upzJd F_1^3 =  \left(  \p_{\olz_1}\phi_1 + \p_{\olz_2}\phi_2 + \p_{\olz_3}\phi_3 + \p_{\olz_4}\phi_4) \right) \gfd_1 \gfd_2 \gfd_3 \gfd_4 I $\\[2mm]
which splits into\\[2mm]
$(\upzJd)_+ F_1^3 =  0$\\[1mm]
$(\upzJd)_- F_1^3 = \left(  \p_{\olz_1}\phi_1 + \p_{\olz_2}\phi_2 + \p_{\olz_3}\phi_3 + \p_{\olz_4}\phi_4) \right) \gfd_1 \gfd_2 \gfd_3 \gfd_4 I$

\end{itemize}
It follows that $F_1^3$ satisfies the following system of equations, which is, as expected, the same system as in case B.
\begin{itemize}

\item $\upz F_1^3 = 0 \Longleftrightarrow \p_{z_1} \phi_2 - \p_{z_2} \phi_1 + \p_{z_3} \phi_4 - \p_{z_4} \phi_3 = 0$
 
\item $ \upzd F_1^3 = 0 \Longleftrightarrow \p_{\olz_2}\phi_4 + \p_{\olz_3}\phi_1 = 0 \ , \  - \p_{\olz_2}\phi_3 +\p_{\olz_4} \phi_1 = 0 \ , \  \p_{\olz_1}\phi_4 - \p_{\olz_3}\phi_2 = 0 \ , \  \p_{\olz_4}\phi_2 + \p_{\olz_1}\phi_3 = 0 \ , \
 \p_{\olz_1}\phi_1 + \p_{\olz_2}\phi_2 = 0 \ , \   \p_{\olz_3}\phi_3 + \p_{\olz_4}\phi_4 = 0$ 
 
\item $ \upzJ F_1^3 = 0 \Longleftrightarrow \p_{z_1}\phi_3 - \p_{z_3}\phi_1 = 0 \ , \ \p_{z_2}\phi_4 - \p_{z_4}\phi_2= 0 \ , \ \p_{z_2}\phi_3 - \p_{z_3}\phi_2 = 0 \ , \ \p_{z_2}\phi_1 - \p_{z_1}\phi_2 = 0 \ , \  \p_{z_1}\phi_4 - \p_{z_4}\phi_1 = 0 \ , \
 \p_{z_4}\phi_3 - \p_{z_3}\phi_4 = 0$
 
\item $\upzJd F_1^3 =  0 \Longleftrightarrow   \p_{\olz_1} \phi_2 + \p_{\olz_2} \phi_1 + \p_{\olz_3} \phi_4 + \p_{\olz_4} \phi_3 = 0$
 
\end{itemize}
We see that $\upzJ F_1^3 = 0 $ implies $\upz F_1^3 = 0 $, but not conversely, and that  $\upzd F_1^3 = 0 $ implies $\upzJd F_1^3 = 0 $, but also not conversely.
We also see that to express q--monogenicity four operators are necessary: $(\upzd)_+$,  $(\upzJ)_+$, either $(\upz)_-$ or $(\upzJ)_-$ which both lead to the same system of equations, and either $(\upzd)_-$ or $(\upzJd)_-$ which also both lead to the same system of equations.\\[-2mm]

\noindent
{\bf Case F:} The function $F_0^4 : \mR^8 \longrightarrow \mS_0^4$ has the form 
$$
F_0^4 = \phi \, \gfd_1 \gfd_2 \gfd_3 \gfd_4 I
$$
There holds
\begin{itemize}
\item $\upzd F_0^4 = \gfd_2 \gfd_3 \gfd_4 \, \p_{\olz_1}\phi -  \gfd_1 \gfd_3 \gfd_4 \, \p_{\olz_2}\phi +  \gfd_1 \gfd_2 \gfd_4 \, \p_{\olz_3}\phi -  \gfd_1 \gfd_2 \gfd_3 \, \p_{\olz_4}\phi$\\[2mm]
which splits into\\[2mm]
$(\upzd)_+ F_0^4 =  \gfd_2 \gfd_3 \gfd_4 \, \p_{\olz_1}\phi -  \gfd_1 \gfd_3 \gfd_4 \, \p_{\olz_2}\phi +  \gfd_1 \gfd_2 \gfd_4 \, \p_{\olz_3}\phi -  \gfd_1 \gfd_2 \gfd_3 \, \p_{\olz_4}\phi$\\[1mm]
$(\upzd)_- F_0^4 = 0$
\end{itemize}
and
\begin{itemize}
\item $\upzJ F_0^4 =  \gfd_2 \gfd_3 \gfd_4 \, \p_{z_2}\phi + \gfd_1 \gfd_3 \gfd_4 \, \p_{z_1}\phi + \gfd_1 \gfd_2 \gfd_4 \, \p_{z_3}\phi + \gfd_1 \gfd_2 \gfd_3 \, \p_{z_4}\phi$\\[2mm]
which splits into\\[2mm]
$(\upzJ)_+ F_0^4 = \gfd_2 \gfd_3 \gfd_4 \, \p_{z_2}\phi + \gfd_1 \gfd_3 \gfd_4 \, \p_{z_1}\phi + \gfd_1 \gfd_2 \gfd_4 \, \p_{z_3}\phi + \gfd_1 \gfd_2 \gfd_3 \, \p_{z_4}\phi$\\[1mm]
$(\upzJ)_- F_0^4 = 0$
\end{itemize}
while trivially
\begin{itemize}
\item $\upzJd F_0^4 = 0 = \upz F_0^4$ 
\end{itemize}
It follows that
\begin{itemize}
\item $\upzd F_0^4 = 0 \Longleftrightarrow \p_{\olz_1}\phi = 0 ,   \p_{\olz_2}\phi = 0 , \p_{\olz_3}\phi = 0 , \p_{\olz_4}\phi = 0$
\item $\upzJ F_0^4 =  0 \Longleftrightarrow  \p_{z_2}\phi = 0 ,  \p_{z_1}\phi = 0 , \p_{z_3}\phi = 0 , \p_{z_4}\phi = 0$
\end{itemize}
in other words: $F_0^4$ is q--monogenic if and only if $\phi$ is a constant function. Also this observation is valid regardless the dimension.
We also see that for expressing the q--monogenicity of $F_0^4$ the two operators $(\upzd)_+$ and $(\upzJ)_+$ need to be used (right upper vertex of the triangle).


\section*{Acknowledgements}

The authors gratefully acknowledge the support by the Czech grant GA CR G201/12/028.



\end{document}